\documentclass{mcom-l}
\usepackage{stmaryrd}

\usepackage{cite}
\usepackage{amsmath}
\usepackage{graphicx,wrapfig}
\usepackage{amssymb,amsgen,color}
\usepackage{amsfonts}
\usepackage{hyperref}
\usepackage{multirow}
\usepackage{epsfig,epstopdf,color,bm}
\usepackage{booktabs}
\usepackage{verbatim}
\usepackage{multicol}
\usepackage{makecell}

\input psfig.sty

\allowdisplaybreaks

\newcommand{\black}{\color{black}}

\renewcommand{\theequation}{\arabic{section}.\arabic{equation}}
\renewcommand{\theequation}{\thesection.\arabic{equation}}
\renewcommand {\thefigure} {\thesection.\arabic{figure}}
\numberwithin{equation}{section}

\newcommand{\bx}{\boldsymbol{x} }

\newcommand{\Og}{\Omega}
\newcommand{\sech}{{\mathrm{sech}} }
\newcommand{\fl}[2]{\frac{#1}{#2}}
\newcommand{\be}{\begin{equation}}
\newcommand{\ee}{\end{equation}}
\newcommand{\nn}{\nonumber}
\newtheorem{theorem}{Theorem}[section]
\newtheorem{lemma}[theorem]{Lemma}
\newcommand{\tabincell}[2]{\begin{tabular}{@{}#1@{}}#2\end{tabular}}

\theoremstyle{definition}

\newtheorem{proposition}[theorem]{Proposition}

\theoremstyle{remark}
\newtheorem{remark}[theorem]{Remark}

\usepackage[top=4cm,bottom=4.0cm,left=3.5cm,right=3.5cm ]{geometry}

\begin{document}

\title[Uniform error bounds of TSFP for NKGE]{Uniform error bounds of  time-splitting spectral methods for the
long-time dynamics of the nonlinear Klein--Gordon equation with weak nonlinearity}
\author{Weizhu Bao}
\address{Department of Mathematics,
National University of Singapore, Singapore 119076}
\email{matbaowz@nus.edu.sg}
\thanks{This work was partially supported by the Ministry of Education of Singapore grant R-146-000-290-114 (W. Bao \& Y. Feng) and Alexander von Humboldt Foundation (C. Su).}

\author{Yue Feng}
\address{(Corresponding author) Department of Mathematics,
National University of Singapore, Singapore 119076}
\email{fengyue@u.nus.edu}

\author{Chunmei Su}
\address{Yau Mathematical Sciences Center, Tsinghua University, 100084 Beijing, China}
\email{sucm@tsinghua.edu.cn}

\date{}

\begin{abstract}
We establish uniform error bounds of time-splitting Fourier pseudospectral (TSFP) methods for the nonlinear Klein--Gordon equation (NKGE) with weak power-type nonlinearity and $O(1)$ initial data, while the nonlinearity strength is characterized by $\varepsilon^{p}$ with a constant $p \in \mathbb{N}^+$ and a dimensionless parameter $\varepsilon \in (0, 1]$, for the long-time dynamics up to the time at $O(\varepsilon^{-\beta})$ with $0 \leq \beta \leq p$. In fact, when $0 < \varepsilon \ll 1$, the problem is equivalent to the long-time dynamics of NKGE with small initial data and $O(1)$ nonlinearity strength, while
the amplitude of the initial data (and the solution) is at $O(\varepsilon)$. By reformulating the NKGE into a relativistic nonlinear Schr\"{o}dinger equation, we adapt the TSFP method to discretize it numerically.
By using the method of mathematical induction to bound the numerical solution, we prove uniform error bounds at  $O(h^{m}+\varepsilon^{p-\beta}\tau^2)$ of the TSFP method with $h$ mesh size, $\tau$  time step and $m\ge2$ depending on the regularity of the solution. The error bounds are  uniformly accurate for the long-time simulation up to the time at $O(\varepsilon^{-\beta})$ and uniformly valid for $\varepsilon\in(0,1]$. Especially, the error bounds are uniformly at the second order rate for the large time step $\tau = O(\varepsilon^{-(p-\beta)/2})$ in the parameter regime $0\le\beta <p$. Numerical results
are reported to confirm our error bounds in the long-time regime. Finally, the TSFP method and its error bounds are extended to a highly oscillatory complex NKGE which propagates waves with wavelength at $O(1)$ in space and $O(\varepsilon^{\beta})$ in time and wave velocity at $O(\varepsilon^{-\beta})$.
\end{abstract}

\maketitle

\section{Introduction}
The nonlinear Klein--Gordon equation (NKGE) is widely used to model nonlinear phenomena in many fields of science and engineering. It plays a fundamental role in quantum electrodynamics, particle and/or plasma physics to describe the motion of spinless particles within the framework of quantum mechanics and Einstein's special relativity \cite{FV,LP,SJJ,SV}. The NKGE with power-type nonlinearity has attracted much attention in investigating the dislocation of crystals, nonlinear optics and quantum field theory \cite{AW,LZ}. In particular,  the NKGE with cubic nonlinearity is called $\varphi^4$ model to describe the relativistic Bose gas, the dynamics of Copper pairs in superconductors as well as displacive and order-disorder transitions in solids \cite{FS2018,KSV}; and the sine-Gordon and sinh-Gordon equations arise in the propagation of fluxons in Josephon junctions between two superconductors \cite{AW}.

In this paper, we consider the following NKGE with power-type nonlinearity on the unit torus $\mathbb{T}^d$ ($d=1,2,3$) as
\be\label{eq:WNE}
\left\{
\begin{aligned}
&\partial_{tt}u(\bx, t)-\Delta u({\bx}, t)+ u({\bx}, t)+\varepsilon^{p}u^{p+1}({\bx}, t)=0,\quad{\bx} \in \mathbb{T}^d,\quad t > 0,\\
&u({\bx}, 0)=u_0({\bx})=O(1),\quad\partial_t u({\bx},0)=u_1({\bx})=O(1),\quad{\bx} \in \mathbb{T}^d.
\end{aligned}\right.
\ee
Here, $t$ is time, $\bx$ is the spatial coordinate, $u:= u(\bx, t)$ is a real-valued scalar field, $p \in \mathbb{N}^+$ is the exponent of the power-type nonlinearity, $\varepsilon\in (0, 1]$ is a dimensionless parameter used to characterize the nonlinearity strength, and the initial datum $u_0({\bx})$ and $u_1({\bx})$ are two given real-valued functions which are independent of the parameter $\varepsilon$. Thus formally,
the amplitude of the solution $u$ is at $O(1)$, the wavelength in space and time is also at $O(1)$, and the wave velocity is at $O(1)$ too. In addition,  if $u(\cdot, t) \in H^1(\mathbb{T}^d)$ and $\partial_t u(\cdot, t) \in L^2(\mathbb{T}^d)$, the NKGE \eqref{eq:WNE} is time symmetric or time reversible and conserves the energy \cite{BD,BFY,DXZ} as
\begin{align}
E(t) &:=E(u(\cdot,t))= \int_{\mathbb{T}^d} \left[ |\partial_t u (\bx, t)|^2 + |\nabla u(\bx, t)|^2 + |u(\bx, t)|^2 +\frac{2\varepsilon^{p}}{p+2} u(\bx, t)^{p+2}  \right] d\bx\nn\\
&\equiv \int_{\mathbb{T}^d} \left[ |u_1(\bx)|^2 + |\nabla u_0(\bx)|^2 + |u_0(\bx)|^2 +\frac{2\varepsilon^{p}}{p+2} u_0(\bx)^{p+2}  \right] d\bx \nn \\
&= E(0) = O(1), \qquad t \geq 0. \label{eq:Energy_u}
\end{align}
In fact, when $0 < \varepsilon \ll 1$, by introducing $w(\bx, t)=\varepsilon u(\bx, t)$, we can reformulate the NKGE \eqref{eq:WNE} with weak nonlinearity and $O(1)$ initial data  into the following NKGE with small initial data and $O(1)$ nonlinearity strength:
\begin{equation}\label{eq:SIE}
\left\{
\begin{aligned}
&\partial_{tt} w({\bx}, t)-\Delta w({\bx}, t)+ w({\bx}, t)+w^{p+1}({\bx}, t) = 0, \quad \bx \in \mathbb{T}^d,\quad t > 0, \\
&w({\bx}, 0) = \varepsilon u_0({\bx})=O(\varepsilon),\quad \partial_t w({\bx}, 0) = \varepsilon u_1({\bx})=O(\varepsilon),\quad {\bx} \in \mathbb{T}^d.
\end{aligned}\right.
\end{equation}
Noticing that the amplitude of the initial data in
\eqref{eq:SIE} is at $O(\varepsilon)$, formally we can get the amplitude of the solution $w$ of \eqref{eq:SIE} is also at $O(\varepsilon)$. Of course,
the wavelength of \eqref{eq:SIE} in space and time  is at $O(1)$, and the wave velocity of \eqref{eq:SIE} is at $O(1)$.
Similarly, the NKGE \eqref{eq:SIE} is time symmetric or time reversible and conserves the energy \cite{BD,BFY,DXZ} as
\begin{align*}
\widetilde{E}(t) &:= \widetilde{E}(w(\cdot,t))= \int_{\mathbb{T}^d} \big[ |\partial_t w (\bx, t)|^2 + |\nabla w(\bx, t)|^2 + |w(\bx, t)|^2 +\frac{2}{p+2} w(\bx, t)^{p+2}  \big] d \bx\\
&=\int_{\mathbb{T}^d} \left[ |\varepsilon u_1(\bx)|^2 + |\varepsilon\nabla u_0(\bx)|^2 + |\varepsilon u_0(\bx)|^2 +\frac{2\varepsilon^{p+2}}{p+2}  u_0(\bx)^{p+2}  \right]  d\bx\\
&=\varepsilon^2 \int_{\mathbb{T}^d} \left[ | u_1(\bx)|^2 + |\nabla u_0(\bx)|^2 + | u_0(\bx)|^2 +\frac{2\varepsilon^{p}}{p+2} u_0(\bx)^{p+2}  \right]  d\bx\\
&=\varepsilon^2 E(0)= O(\varepsilon^2), \qquad t\ge0.
\end{align*}
Thus, the long-time dynamics of the NKGE \eqref {eq:SIE} with small initial data and $O(1)$ nonlinearity strength is equivalent to the long-time dynamics of the NKGE \eqref{eq:WNE} with weak nonlinearity and $O(1)$ initial data. In both cases,
the solutions propagate waves with wavelength in space and time at
$O(1)$ and the wave velocity at $O(1)$.

There are two different dynamical problems related to the time evolution of the NKGE
\eqref{eq:WNE} (or \eqref{eq:SIE}): (i) when $\varepsilon =\varepsilon_0$ (e.g., $\varepsilon=1$) fixed, i.e., in the standard nonlinearity strength regime, to study the finite time dynamics of \eqref{eq:WNE} (or \eqref{eq:SIE}) for $t\in[0,T]$ with $T=O(1)$; and (ii) when $0<\varepsilon\ll1$, i.e., in the weak nonlinearity strength regime, to study the long-time dynamics of \eqref{eq:WNE} (or \eqref{eq:SIE}) for $t\in[0,T_\varepsilon]$ with $T_\varepsilon=O(\varepsilon^{-{p}})$. Extensive mathematical and numerical studies have been done in the literature for the finite time dynamics of \eqref{eq:WNE} with $\varepsilon=1$, i.e., in the standard nonlinearity strength regime. Along the analytical front, for the existence of global classical solutions, approximate and almost periodic solutions as well as asymptotic behavior of the solution of \eqref{eq:WNE} with $\varepsilon=1$,
we refer to \cite{BJ,BV,CE,KT,K,ONO,VW} and references therein. For the numerical aspects, different numerical methods have been presented and analyzed in the literature, such as finite difference time domain (FDTD) methods,  spectral methods, etc. For details, we refer to \cite{BCZ,BD,CCLM,DXZ,DB,FS} and references therein. Recently, there are several analytical studies for the long-time dynamics of \eqref{eq:WNE} in the weak nonlinearity strength regime (or \eqref{eq:SIE} with small initial data), i.e., $0<\varepsilon\ll1$ \cite{KT,LH}. According to the analytical results, the life-span of a smooth solution to the NKGE \eqref{eq:WNE} (or \eqref{eq:SIE}) is at least up to the time at $O(\varepsilon^{-p})$ \cite{D,D2,DS,FZ,KT,K}.

However, to the best of our knowledge, there are very few numerical analytical results on error bounds of the numerical methods for the long-time dynamics of \eqref{eq:WNE} in the literature, especially the error bounds which are valid up to the time at $T_\varepsilon =O(\varepsilon^{-p})$ and how the error bounds depend explicitly  on the mesh size $h$ and time step $\tau$ as well as the small parameter $\varepsilon\in (0,1]$. We notice that some numerical analysis results on the long-time near-conservation (or approximate preservation) of energy, momentum and harmonic actions have been established for some semi-discretizations or full discretizations of the NKGE \eqref{eq:SIE} with small initial data  via the technique of modulated Fourier expansions \cite{CHL,CHL2,HL}, however, no error estimate of the numerical solution itself has been given in the literature. Recently, for the NKGE \eqref{eq:WNE} with cubic nonlinearity (i.e., $p=2$), error estimates of four different FDTD methods were established for the long-time dynamics of the NKGE \eqref{eq:WNE} up to the long-time at $O(\varepsilon^{-\beta})$
with $0 < \beta \leq 2$ \cite{BFY, Feng}. Specifically, in order to obtain `correct' numerical approximations of the NKGE \eqref{eq:WNE} (or \eqref{eq:SIE}) up to the long-time at $O(\varepsilon^{-\beta})$ with $0 < \beta \leq 2$, the $\varepsilon$-scalability (or meshing strategy) of the FDTD methods should be
\be\label{msfdtd}
h=O(\varepsilon^{\beta/2})\quad \mbox{and}\quad \tau=O(\varepsilon^{\beta/2}),
\ee
which immediately suggests that the FDTD methods are  {\bf under-resolution} in both space and time with respect to $\varepsilon\in(0,1]$ in terms of the resolution capacity of the
Shannon's information theory \cite{Lan,Shan1} -- to resolve a wave one needs a few points per wave --
since the wavelength of the solution of the NKGE \eqref{eq:WNE} (or \eqref{eq:SIE}) in space and time is at $O(1)$, while the mesh size $h$ and time step $\tau$ have to be taken at $O(\varepsilon^{\beta/2})$ which is much smaller than $O(1)$! In fact, the FDTD methods can also be regarded as {\bf over-sampling} methods in the sense that the number of points needed per wave in space and time have to be taken as $O(\varepsilon^{-\beta/2})$
which is much larger than $O(1)$! To improve this, a Gautschi-type exponential wave integrator Fourier pseudospectral (EWI-FP) method was proposed and analyzed in \cite{FY}, where a uniform error bound was established at $O(h^{m}+\varepsilon^{2-\beta}\tau^2)$ under a stability condition $\tau\lesssim h$, while $m\ge2$ depending on the regularity of the solution, for the long-time dynamics up to the time at $O(\varepsilon^{-\beta})$ with $0<\beta\le 2$.

As we know, the time-splitting Fourier pseudospectral (TSFP) method has been
widely used to numerically solve dispersive partial differential equations (PDEs) \cite{BC,BCJT,BCJY,BS,DXZ,HJMSZ,LU,SW}. In many cases, the TSFP method demonstrates much better spatial/temporal resolution than the FDTD methods, especially when they are used for integrating highly oscillatory PDEs, such as for the Schr\"{o}dinger/nonlinear Schr\"{o}dinger equation in the semiclassical regime \cite{BJM,CG}, for the NKGE in the nonrelativistic regime \cite{DXZ}, for the Zakharov system in the subsonic limit regime \cite{BSu}, for the Dirac/nonlinear Dirac equation in the nonrelativistic regime \cite{BCJT,BCJY}, etc. The main aim of this paper is to adapt the TSFP method for discretizing the NKGE \eqref{eq:WNE} and establish its error bound for the long-time dynamics up to the time at $O(\varepsilon^{-\beta})$ with $0 \leq \beta \leq p$. In order to do so, we first reformulate the NKGE \eqref{eq:WNE} into a relativistic nonlinear Schr\"{o}dinger equation (NLSE) and then apply the TSFP method to discretize it numerically.  By employing the method of mathematical induction to bound the numerical solution, we establish an error bound at $O(h^{m}+\varepsilon^{p-\beta}\tau^2)$ without any stability condition, while $m\ge2$ depends on the regularity of the solution, for the long-time dynamics up to the time at $O(\varepsilon^{-\beta})$ with $0\le \beta\le p$. The error bound immediately indicates that the TSFP method is uniformly accurate for the long-time simulation up to the time at $O(\varepsilon^{-p})$ and is uniformly valid for $\varepsilon \in (0,1]$. Thus, the TSFP method is an {\bf optimal resolution} method for the long-time dynamics of the NKGE \eqref{eq:WNE} up to the time at $O(\varepsilon^{-p})$. Compared to the EWI-FP method in \cite{FY}, the TSFP method is superior on several aspects: (i) the strict stability condition $\tau\lesssim h$ is removed, (ii) the error bounds are uniformly second order accurate for the large time step $\tau = O(\varepsilon^{-(p-\beta)/2})$ when $0\le \beta <p$, and (iii) we observe numerically the TSFP method has an improved convergence when $0<\varepsilon \ll 1$, which is not valid for the EWI-FP method (cf. Sect. 4).

The rest of the paper is organized as follows. In Sect. 2, we first reformulate the NKGE \eqref{eq:WNE} into a relativistic NLSE and then present the TSFP method to discretize it numerically. In Sect. 3, we establish uniform error bounds of the TSFP
method for the long-time dynamics of the NKGE \eqref{eq:WNE} up to time
at $O(\varepsilon^{-\beta})$ with $0 \leq \beta \leq p$. Numerical results are reported in Sect. 4 to confirm the error estimates. Extension to a highly oscillatory complex  NKGE  in the whole space is presented in  Sect. 5. Finally, some conclusions are drawn in Sect. 6. Throughout this paper, $C$ represents a generic constant which is independent of the discretization parameters $h$ and $\tau$ as well as the nonlinearity strength parameter $\varepsilon\in(0,1]$. We adopt the notation $A \lesssim B$ to represent that there exists a generic constant $C>0$ such that $|A| \leq C B$, while $C$ is independent of $h$ and $\tau$ as well as $\varepsilon$.

\section{A time-splitting Fourier pseudospectral (TSFP) method}
In this section, we first reformulate the NKGE \eqref{eq:WNE} into
a relativistic NLSE and then adopt
the TSFP method \cite{BC,BS,DXZ,HJMSZ,LU,YRS}  to discretize it numerically.

\subsection{A relativistic nonlinear Schr\"{o}dinger equation (NLSE) }

For simplicity of notations, we only illustrate the approach in one dimension (1D) and all the notations and results can be easily generalized to higher dimensions with minor modifications. In 1D, the NKGE \eqref{eq:WNE} with periodic boundary condition collapses to
\be\label{eq:21}
\left\{
\begin{aligned}
&\partial_{tt} u(x, t) - \partial_{xx} u (x, t)+  u(x, t) + \varepsilon^{p} u^{p+1} (x, t)= 0,\quad x \in \Omega = (a, b),\quad t > 0, \\
&u(a,t)=u(b,t), \qquad \partial_x u(a,t)=\partial_x u(b,t), \qquad
t\ge0,\\
&u(x, 0) =u_0(x), \qquad \partial_t u(x, 0) =u_1(x) , \qquad x \in \overline{\Omega} = [a, b].
\end{aligned}\right.
\ee

 For an integer $m\ge 0$,
$\Omega=(a, b)$, we denote by $H^m(\Omega)$  the standard Sobolev space with norm
\be\label{sn}
\|z\|_m^2=\sum\limits_{l \in \mathbb{Z}} (1+|\mu_l|^2)^m|\widehat{z}_l|^2,\quad \mathrm{for}\quad z(x)=\sum\limits_{l\in \mathbb{Z}} \widehat{z}_l e^{i\mu_l(x-a)},\quad \mu_l=\fl{2\pi l}{b-a},
\ee
where $\widehat{z}_l (l\in \mathbb{Z})$ are the Fourier transform coefficients  of the function $z(x)$ \cite{BCJY,BCZ}.
For $m=0$, the space is exactly $L^2(\Og)$ and the corresponding norm is denoted as $\|\cdot\|$.
Furthermore, we denote by $H_{\rm per}^m(\Omega)$ the subspace of $H^m(\Omega)$ which consists of functions with derivatives of order up to $m-1$ being $(b-a)$-periodic.
We see that the space $H^m(\Omega)$ with fractional $m$ is also well-defined which consists of functions with finite norm $\|\cdot\|_m$ \cite{ST}.

Define the operator
\begin{equation}
\langle \nabla \rangle=\sqrt{1-\Delta},
\end{equation}
through its action in the Fourier space by \cite{FS,SZ}:
\[
\langle \nabla \rangle z(x)=\sum\limits_{l\in\mathbb{Z}}\sqrt{1+|\mu_l|^2}\widehat{z}_l e^{i\mu_l(x-a)}, \quad \mathrm{for}\quad z(x)=\sum\limits_{l\in\mathbb{Z}} \widehat{z}_l e^{i\mu_l(x-a)},\quad x\in[a,b].\]
Then we can rewrite the NKGE \eqref{eq:21} as
\begin{equation}
\partial_{tt} u(x, t) + \langle\nabla \rangle^2 u(x, t) + \varepsilon^{p} u^{p+1}(x, t) = 0,\quad x \in\Omega,\quad t>0.
\end{equation}
In addition, we introduce the operator $ \langle \nabla \rangle^{-1}$  as
\[\langle \nabla \rangle^{-1} z(x)=\sum\limits_{l\in\mathbb{Z}}\frac{\widehat{z}_l}{\sqrt{1+|\mu_l|^2}} e^{i\mu_l(x-a)},\qquad x\in\overline{\Omega}.\]
It is obvious that
\[\| \langle \nabla \rangle^{-1} z\|_s =\|z\|_{s-1}\le \|z\|_s.\]

Denote $v(x, t) = \partial_t u(x, t)$ and set
\begin{equation}
\psi(x, t) = u(x, t) - i\langle\nabla\rangle^{-1}v(x, t),
\quad x\in[a,b], \quad t\ge0.
\label{eq:psi}
\end{equation}
By a short calculation, we can reformulate the NKGE \eqref{eq:21} into a relativistic NLSE
in  $\psi:=\psi(x,t)$ as
\be\label{eq:NLS}
\left\{
\begin{aligned}
&i\partial_t \psi(x, t) + \langle\nabla \rangle \psi(x, t) + \ \varepsilon^{p}\langle\nabla \rangle^{-1} f\Big(\frac{1}{2} \left(\psi + \overline{\psi}\right)\Big)(x, t) = 0,\,\,\,\,\, x \in \Omega, \,\,\,\,\, t > 0,\\
&\psi(a,t)=\psi(b,t), \qquad \partial_x \psi(a,t)=\partial_x \psi(b,t), \quad
t\ge0,\\
&\psi(x,0)=\psi_0(x):=u_0(x)-i\langle\nabla \rangle^{-1}u_1(x), \quad x\in[a,b],
\end{aligned}\right.
\ee
where $f(z)=z^{p+1}$ and $\overline{\psi}$ denotes the complex conjugate of $\psi$. Noticing \eqref{eq:psi}, we can recover the solution of the NKGE \eqref{eq:21} by
\be\label{uxtpsit}
u(x,t)=\frac{1}{2}\left(\psi(x,t)+\overline{\psi}(x,t)\right), \qquad
v(x, t)=\frac{i}{2}\langle\nabla \rangle\left(\psi(x,t)-\overline{\psi}(x,t)\right).
\ee

We remark here that the NKGE \eqref{eq:21} can also be reformulated
as the following first-order (in time) PDEs:
\be\label{eq:split278}
\left\{
\begin{aligned}
&\partial_t u(x, t) - v(x, t) = 0, \quad x \in (a,b), \quad t > 0,\\
&\partial_t v(x, t) - \partial _{xx} u(x, t) + u(x, t)+ \varepsilon^{p} u^{p+1}(x, t)= 0, \quad x \in (a,b), \quad t > 0, \\
&u(a,t)=u(b,t), \qquad \partial_x u(a,t)=\partial_x u(b,t), \qquad
t\ge0,\\
&u(x,0)=u_0(x),\quad v(x,0) = u_1(x), \qquad x\in[a,b].
\end{aligned}\right.
\ee

\subsection{Semi-discretization by using the second-order time-splitting}

In order to discretize the NKGE \eqref{eq:21} in time by a time-splitting method,  we first discretize the relativistic NLSE  \eqref{eq:NLS}
by a time-splitting method and then recover the solution of \eqref{eq:21}
via \eqref{uxtpsit}. In fact, the relativistic NLSE  \eqref{eq:NLS}
can be decomposed into the following two subproblems via  the
time-splitting technique \cite{LU,SZ}
\be\label{eq:NLS_s2}
\left\{
\begin{aligned}
&i\partial_t\psi(x,t)+\langle\nabla\rangle \psi(x,t)=0,\quad x \in (a,b), \quad t > 0,\\
&\psi(a,t)=\psi(b,t), \qquad \partial_x \psi(a,t)=\partial_x \psi(b,t), \quad
t\ge0,\\
&\psi(x,0)=\psi_0(x), \qquad x\in[a,b],
\end{aligned}\right.
\ee
and
\be\label{eq:NLS_s1}
\left\{
\begin{aligned}
&i\partial_t \psi(x,t) +\varepsilon^{p}\langle\nabla \rangle^{-1} f\Big(\frac{1}{2} (\psi+ \overline{\psi})\Big)(x,t)=0,\quad x \in (a, b), \quad t > 0,\\
&\psi(x,0)=\psi_0(x),\qquad x\in[a,b].
\end{aligned}\right.
\ee
The linear equation \eqref{eq:NLS_s2} can be solved {\sl exactly} in phase space and the associated evolution operator is given by
\be\label{lex1}
\psi(\cdot, t)=\varphi^{t}_T(\psi_0):=e^{it\langle\nabla\rangle}\psi_0,\quad t\ge 0,
\ee
which satisfies the isometry relation
\[\|\varphi^{t}_T(v_0)\|_s=\|v_0\|_s,\quad s\ge 0,\quad t\in\mathbb{R}.\]
Recalling that the nonlinear part of \eqref{eq:NLS_s1} is real, this implies that $\partial_t \left(\psi+\overline \psi \right)(x,t)=0$ for any fixed
$x\in[a,b]$. Thus $\psi+\overline \psi$ is invariant in time, i.e.,
\be\label{conspsi}
\left(\psi+\overline \psi \right)(x,t)\equiv \left(\psi+\overline \psi \right)(x,0)=\psi_0(x)+\overline{\psi_0}(x),
\qquad t\ge0, \quad a\le x\le b.
\ee
Plugging \eqref{conspsi} into \eqref{eq:NLS_s1}, we get
\be\label{eq:NLS_s21}
\left\{
\begin{aligned}
&i\partial_t \psi(x,t) +\varepsilon^{p}\langle\nabla \rangle^{-1} f\Big(\frac{1}{2} (\psi_0+ \overline{\psi_0})\Big)(x)=0,\quad x \in [a,b], \quad t > 0,\\
&\psi(x,0)=\psi_0(x),\qquad x\in[a,b].
\end{aligned}\right.
\ee
Thus \eqref{eq:NLS_s21} (or \eqref{eq:NLS_s1}) can be integrated
{\sl exactly} in time as:
\be\label{nlex1}
\psi(x,t)=\varphi^{t}_V (\psi_0) := \psi_0(x) + \varepsilon^{p}t \, F(\psi_0(x)),\quad t\ge0,\ee
where the operator $F$ is defined by
\be\label{Fd}
F(\phi)=i\langle\nabla\rangle^{-1} G(\phi),\qquad G(\phi)=f\Big(\frac{1}{2}(\phi + \overline{\phi})\Big).
\ee

Let $\tau>0$ be the time step and define $t_n=n\tau$ for $n=0, 1,\ldots$. Denote $\psi^{[n]}:=\psi^{[n]}(x)$ by the approximation
of $\psi(x,t_n)$ for $n\ge0$, then a second-order semi-discretization
of the relativistic NLSE \eqref{eq:NLS} via the Strang splitting \cite{LU} can be given as:
\be\label{Strang}
\psi^{[n+1]}=\mathcal{S}_{\tau}(\psi^{[n]})=\varphi^{\tau/2}_T \circ\varphi^{\tau}_V \circ\varphi^{\tau/2}_T(\psi^{[n]}),\qquad n=0,1,2,\ldots,
\ee
with $\psi^{[0]}=\psi_0=u_0-i\langle\nabla\rangle^{-1}u_1$.
Noticing \eqref{uxtpsit} and \eqref{Strang}, we can get a second-order semi-discretization of the NKGE \eqref{eq:21}:
\be\label{uxtpsita}
u^{[n]}=\frac{1}{2}\left(\psi^{[n]}+\overline{\psi^{[n]}}\right), \quad
v^{[n]} = \frac{i}{2}\langle\nabla \rangle\left(\psi^{[n]}-\overline{\psi^{[n]}}\right),\quad n=0, 1,\ldots,
\ee
where $u^{[n]}:=u^{[n]}(x)$ and $v^{[n]}:= v^{[n]}(x)$
are the approximations of $u(x,t_n)$ and $\partial_t u(x,t_n)$ ($n=0,1,2,\ldots$), respectively.

We remark here that another way to discretize the NKGE \eqref{eq:21} by a time-splitting method, which is exactly the same discretization as the one presented above, is to discretize the NKGE \eqref{eq:split278} by a
time-splitting method. In fact, the NKGE \eqref{eq:split278}
can be decomposed into the following two subproblems via  the
time-splitting technique \cite{DXZ}
\be\label{eq:split2}
\left\{
\begin{aligned}
&\partial_t u(x, t) - v(x, t) = 0, \\
&\partial_t v(x, t) - \partial _{xx} u(x, t) + u(x, t) = 0, \quad x \in (a,b), \quad t > 0,\\
&u(a,t)=u(b,t), \qquad \partial_x u(a,t)=\partial_x u(b,t), \qquad
t\ge0,\\
&u(x,0)=u_0(x),\quad v(x,0) = u_1(x), \qquad x\in[a,b],
\end{aligned}\right.
\ee
and
\be\label{eq:split1}
\left\{
\begin{aligned}
&\partial_t u(x, t) = 0, \\
&\partial_t v(x, t)  +\varepsilon^{p} u^{p+1}(x, t) = 0, \quad x \in [a,b], \quad t > 0,\\
&u(x,0)=u_0(x),\quad v(x,0) = u_1(x), \qquad x\in[a,b].
\end{aligned}\right.
\ee
Similarly, the linear problem \eqref{eq:split2} can be solved {\sl exactly} in phase space and the associated evolution operator is given by
\be\label{lex}
\left(
\begin{aligned}
&u(\cdot, t)\\
&v(\cdot, t)
\end{aligned}\right)
=\chi^{t}_T
\left(
\begin{aligned}
&u_0\\
&u_1
\end{aligned}\right):=
\left(
\begin{aligned}
&\cos(t\langle \nabla\rangle)u_0+\langle \nabla \rangle^{-1}\sin(t\langle \nabla\rangle)u_1 \\
&-\langle \nabla\rangle\sin(t\langle \nabla\rangle)u_0+\cos(t\langle \nabla\rangle)u_1
\end{aligned}\right),\quad t\ge 0.
\ee
From \eqref{eq:split1}, we obtain immediately that $u(x,t)$ is invariant in time for any fixed $x\in[a,b]$, i.e.,
\be\label{uxt987}
u(x,t)\equiv u(x,0)=u_0(x), \qquad x\in[a,b].
\ee
Plugging \eqref{uxt987} into \eqref{eq:split1}, we get
\be\label{eq:split187}
\left\{
\begin{aligned}
&\partial_t u(x, t) = 0, \\
&\partial_t v(x, t) + \varepsilon^{p}u^{p+1}(x, 0)= 0, \quad x \in [a,b], \quad t > 0,\\
&u(x,0)=u_0(x),\quad  v(x,0) = u_1(x), \quad t\ge0,  \qquad x\in[a,b].
\end{aligned}\right.
\ee
Thus \eqref{eq:split187} (and \eqref{eq:split1}) can be integrated {\sl exactly} in time as:
\be\label{nkgesp345}
\left(
\begin{aligned}
&u(\cdot, t)\\
&v(\cdot, t)
\end{aligned}\right)
=\chi^{t}_V
\left(
\begin{aligned}
&u_0\\
&u_1
\end{aligned}\right):=
\left(
\begin{aligned}
&u_0\\
&u_1- \varepsilon^{p} t u_0^{p+1}
\end{aligned}\right),\quad t\ge 0.
\ee
Let $u^{[n]}:=u^{[n]}(x)$ and $v^{[n]}: = v^{[n]}(x)$
be the approximations of $u(x,t_n)$ and $v(x,t) = \partial_t u(x,t_n)$ ($n=0,1,2,\ldots$), respectively, which are the solutions of the NKGE \eqref{eq:split278} (and \eqref{eq:21}). Then a second-order semi-discretization
of the NKGE \eqref{eq:split278} (and \eqref{eq:21}) via the second-order Strang splitting \cite{DXZ} can be given as:
\be\label{Strang1}
\left(
\begin{aligned}
&u^{[n+1]}\\
&v^{[n+1]}
\end{aligned}\right)
=\mathcal{S}_{\tau}\left(
\begin{aligned}
&u^{[n]}\\
&v^{[n]}
\end{aligned}\right)=\chi^{\tau/2}_T \circ\chi^{\tau}_V \circ\chi^{\tau/2}_T
\left(
\begin{aligned}
&u^{[n]}\\
&v^{[n]}
\end{aligned}\right), \quad n=0, 1, \ldots,
\ee
with $u^{[0]}=u_0$ and $v^{[0]}=u_1$. In fact, it is easy  to verify that \eqref{eq:NLS_s2}, \eqref{eq:NLS_s1}, \eqref{lex1} and \eqref{nlex1} are equivalent to \eqref{eq:split2}, \eqref{eq:split1}, \eqref{lex} and \eqref{nkgesp345}, respectively. Thus it is straightforward to get that
\eqref{Strang} is equivalent to \eqref{Strang1}, and \eqref{uxtpsita} is the same as \eqref{Strang1}.

\begin{remark} Another second-order semi-discretization of the relativistic NLSE \eqref{eq:NLS} can be given as
\be\label{Deuf}
\psi^{[n+1]}=\varphi^{\tau/2}_V \circ\varphi^{\tau}_T \circ\varphi^{\tau/2}_V(\psi^{[n]}), \qquad n=0,1,2,\ldots\;,
\ee
which can immediately generate a semi-discretization of the NKGE \eqref{eq:21}  via \eqref{uxtpsita}. Again, it is easy to check that
this discretization is the same as the discretization of the NKGE
\eqref{eq:split278} (and \eqref{eq:21}) by a similar second-order Strang-type time-splitting as
\be
\left(
\begin{aligned}
&u^{[n+1]}\\
&v^{[n+1]}
\end{aligned}\right)
=\chi^{\tau/2}_V \circ\chi^{\tau}_T \circ\chi^{\tau/2}_V
\left(
\begin{aligned}
&u^{[n]}\\
&v^{[n]}
\end{aligned}\right), \qquad n=0,1,2,\ldots.
\ee
Furthermore, the above second-order time-splitting discretization of the
NKGE \eqref{eq:21} is equivalent to an exponential wave integrator
(EWI) via the trapezoidal quadrature (or Deuflhard-type exponential integrator) for discretizing the NKGE \eqref{eq:21} directly (cf. \cite{DXZ}).
\end{remark}

\begin{remark}
It is straightforward to design higher order semi-discretizations of the NKGE \eqref{eq:21} via the relativistic NLSE \eqref{eq:NLS} by adopting a
higher order time-spitting method \cite{MQ}, e.g., the fourth-order partitioned Runge-Kutta time-splitting method \cite{BS}.
\end{remark}

\subsection{Full-discretization by the Fourier pseudospectral method}
Let $N$ be an even  positive integer and define the spatial mesh size $h=(b-a)/N$, then the grid points are chosen as
\begin{equation}
x_j := a + jh,\quad j \in \mathcal{T}^0_N=\{j~|~j = 0, 1, \ldots, N\}.
\end{equation}
Denote $X_N := \{z= (z_0, z_1, \ldots, z_N)^T \in \mathbb{R}^{N+1} \ | \ z_0 = z_N\}$ with the $l^2$-norm and $l^{\infty}$-norm in $X_N$ given as
\begin{equation}
\|z\|^2_{l^2} = h \sum_{j=0}^{N-1} |z_j|^2,\quad \|z\|_{l^{\infty}} = \max_{0\le j \le N-1} |z_j|, \quad z \in X_N.
\end{equation}
Define $C_{\rm per}(\Omega)=\{z\in C(\overline \Omega) \ |\ z(a)=z(b)\}$ and
\[
Y_N := \text{span}\left\{e^{i\mu_l(x-a)},\quad x \in \overline{\Omega}, \quad l \in \mathcal{T}_N\right\},\quad \mathcal{T}_N = \left\{l ~|~ l = -\frac{N}{2}, -\frac{N}{2}+1, \ldots, \frac{N}{2}-1\right\}.
\]
For any $z(x) \in C_{\rm per}(\Omega)$ and a vector $z\in X_N$, let $P_N: L^2(\Omega) \to Y_N$ be the standard $L^2$-projection operator onto $Y_N$, $I_N : C_{\rm per}(\Omega) \to Y_N$ or $I_N : X_N \to Y_N$ be the trigonometric interpolation operator \cite{ST}, i.e.,

\begin{equation}
(P_N z)(x) = \sum_{l \in \mathcal{T}_N} \widehat{z}_l e^{i\mu_l(x-a)},\qquad (I_N z)(x) = \sum_{l \in \mathcal{T}_N} \widetilde{z}_l e^{i\mu_l(x-a)},\qquad x \in \overline{\Omega},
\end{equation}
where
\begin{equation}
\widehat{z}_l = \frac{1}{b-a}\int^b_a z(x) e^{-i\mu_l (x-a)} dx, \quad \widetilde{z}_l = \frac{1}{N}\sum_{j=0}^{N-1} z_j e^{-i\mu_l (x_j-a)}, \quad l \in \mathcal{T}_N,
\end{equation}
with $z_j$ interpreted as $z(x_j)$ when involved.

Let $\psi_j^n$ be the numerical approximation of $\psi(x_j,t_n)$
for $j\in \mathcal{T}^0_N$ and $n\ge0$ and denote $\psi^n=(\psi_0^n, \psi_1^n,\ldots, \psi_N^n)^T\in \mathbb{C}^{N+1}$ for $n=0,1,\ldots$. Then a time-splitting
Fourier pseudospectral (TSFP) method for discretizing the relativistic NLSE  \eqref{eq:NLS} via \eqref{Strang} with a Fourier pseudospectral discretization in space can be given as
\be\label{psifull}
\begin{split}
&\psi^{(n, 1)}_j=\sum_{l \in \mathcal{T}_N} e^{i\frac{\tau\zeta_l}{2}}\;\widetilde{(\psi^n)}_l\; e^{i\mu_l(x_j-a)},  \\
&\psi^{(n, 2)}_j=\psi^{(n, 1)}_j+\varepsilon^{p}\tau\, F_j^n, \qquad F_j^n=i\sum_{l \in \mathcal{T}_N} \frac{1}{\zeta_l} \widetilde{\left(G(\psi^{(n,1)})\right)}_l\; e^{i\mu_l(x_j-a)}, \\
&\psi^{n+1}_j=\sum_{l \in \mathcal{T}_N} e^{i\frac{\tau\zeta_l}{2}} \; \widetilde{\left(\psi^{(n, 2)}\right)}_l\; e^{i\mu_l(x_j-a)},\quad j \in \mathcal{T}^0_N, \quad n=0,1,\ldots,
\end{split}
\ee
where $\zeta_l=\sqrt{1+\mu_l^2}$ for $l\in \mathcal{T}_N$, $\psi^{(n,k)}=(\psi_0^{(n,k)}, \psi_1^{(n,k)},\ldots$, $\psi_N^{(n,k)})^T\in \mathbb{C}^{N+1}$ for
$k=1$, $2$, $G(\psi^{(n,1)}):=(G(\psi^{(n,1)}_0), G(\psi^{(n,1)}_2), \ldots, G(\psi^{(n,1)}_N))^T\in \mathbb{R}^{N+1}$ and
\[\psi_j^0=u_0(x_j)-i\sum_{l \in \mathcal{T}_N}\frac{\widetilde{(u_1)}_l}{\sqrt{1+|\mu_l|^2}} e^{i\mu_l(x_j-a)}, \qquad j \in \mathcal{T}^0_N. \]

Let $u^n_j$ and $v^n_j$ be the approximations of $u(x_j, t_n)$ and $v(x_j, t_n)$, respectively, for $j\in \mathcal{T}^0_N$ and $n\ge0$,
 and denote $u^n=(u_0^n, u_1^n,\ldots, u_N^n)^T\in \mathbb{R}^{N+1}$
and  $v^n=(v_0^n, v_1^n,\ldots, v_N^n)^T\in \mathbb{R}^{N+1}$  for $n=0,1,\ldots$.
Combining \eqref{psifull} and \eqref{uxtpsita}, we can obtain a
full-discretization of the NKGE \eqref{eq:21} by the TSFP method as
\be\label{ufull}
\begin{split}
&u_j^{n+1}=\frac{1}{2}\left(\psi_j^{n+1}+\overline{\psi_j^{n+1}}\right),\\
&v_j^{n+1} = \frac{i}{2}\sum_{l \in \mathcal{T}_N}\zeta_l\big[\widetilde{(\psi^{n+1})}_l-
\widetilde{(\overline{\psi^{n+1}})}_l\big]\;
e^{i\mu_l (x_j-a)},
\end{split}
\qquad j \in \mathcal{T}^0_N, \quad n\ge 0,
\ee
with
\[u_j^0=u_0(x_j), \qquad v_j^0 = u_1(x_j), \qquad j \in \mathcal{T}^0_N.
\]

Specifically, plugging \eqref{psifull} into \eqref{ufull} or discretizing \eqref{Strang1}
directly in space by the Fourier pseudospectral method, we get a
full-discretization of the NKGE \eqref{eq:21} by the TSFP method (in explicit formulation in the original variable $u$) as
\be\label{ufull2}
\begin{split}
&u^{(n, 1)}_j =\mathcal{L}_u\left(\frac{\tau}{2}, u^n, v^n\right)_j,  \qquad\qquad \,\,v^{(n, 1)}_j = \mathcal{L}_v\left(\frac{\tau}{2}, u^n, v^n\right)_j,  \\
&u^{(n, 2)}_j = u^{(n, 1)}_j,\qquad\qquad\qquad\qquad\quad v^{(n, 2)}_j = v^{(n, 1)}_j- \tau \varepsilon^{p}\big(u^{(n, 1)}_j\big)^{p+1},\\
&u^{n+1}_j = \mathcal{L}_u\left(\frac{\tau}{2}, u^{(n, 2)}, v^{(n, 2)}\right)_j,\qquad v^{n+1}_j = \mathcal{L}_v \left(\frac{\tau}{2}, u^{(n, 2)}, v^{(n, 2)}\right)_j,
\end{split}
\ee
where
\be\label{Lus}
\begin{split}
&\mathcal{L}_u \left(\tau, u, v\right)_j = \sum_{l \in \mathcal{T}_N} \left[\cos(\tau\zeta_l)\widetilde{u}_l + \zeta^{-1}_l \sin(\tau\zeta_l)\widetilde{v}_l \right]e^{i\mu_l(x_j-a)},\\
&\mathcal{L}_v \left(\tau, u, v\right)_j = \sum_{l \in \mathcal{T}_N} \left[-\zeta_l\sin(\tau\zeta_l)\widetilde{u}_l  + \cos(\tau\zeta_l)\widetilde{v}_l \right]e^{i\mu_l(x_j-a)},
\end{split}
\quad j \in \mathcal{T}^0_N.
\ee

The TSFP method \eqref{ufull2} (or \eqref {ufull} with \eqref{psifull}) for the NKGE \eqref{eq:21} is explicit, time symmetric and easy to be extended to higher dimensions. The memory cost of the TSFP method is $O(N)$ and the computational cost per time step is $O(N\ln N)$. In addition, the total
cost for the long-time dynamics up to the time $T_\varepsilon=T_0/\varepsilon^\beta$ ($0\leq \beta \le p$) with fixed $T_0>0$ is
$O\left(\frac{T_\varepsilon\, N\ln N}{\tau}\right)=O\left(\frac{T_0N\ln N}{ \tau\varepsilon^\beta}\right)$.

\section{Uniform error bounds of the TSFP method}
In this section, we establish error bounds of the TSFP method
\eqref{ufull} via \eqref{psifull} (or equivalently \eqref{ufull2}) for the NKGE \eqref{eq:21} up to the time at $O(\varepsilon^{-p})$, which are uniformly valid for $0<\varepsilon\le 1$.

\subsection{Main results}
Motivated by the discussions in \cite{D,FZ,ONO} and references therein, we make the following assumptions on the exact solution $u:=u(x, t)$ of the NKGE \eqref{eq:21} up to the time at $T_{\varepsilon} = T_0/\varepsilon^\beta$ with $\beta\in [0, p]$ and $T_0>0$ fixed:
\[
{\rm(A)}
\begin{split}
&u \in  \  L^\infty\left([0, T_{\varepsilon}]; H^{m+1}_{\rm per}\right), \qquad
\partial_t u\in L^\infty\left([0, T_{\varepsilon}]; H^{m}_{\rm per}\right),\\
&\|u\|_{L^{\infty}\left([0, T_{\varepsilon}]; H^{m+1}_{\rm per}\right)} \lesssim 1,\qquad\,\,\,\, \|\partial_t u\|_{L^{\infty}\left([0, T_{\varepsilon}]; H^{m}_{\rm per} \right)} \lesssim 1,
\end{split}
\]
with $m\ge 1$. Then we can establish the following error bounds of the TSFP method.

\begin{remark}
For the quadratic nonlinearity, i.e., $p=1$, the assumption (A) can be established under the condition on the initial data satisfying
\[
u_0 \in  \  H^{m+1}_{\rm per}, \qquad
u_1 \in H^{m}_{\rm per},\]
if $m > \frac{d}{2}+1$ with $d$ representing the dimension of the torus \cite{D2}. For $p >1$, the regularity of the solution $u(x, t)$ can be preserved and the uniform boundedness in (A) can be established up to the time until $T_\varepsilon=T_0/\varepsilon^{p}$ when $m$ is large enough \cite{D,BFG,FZ}.
\end{remark}

\begin{theorem}\label{thm:eb1}
Let $u^n$ be the numerical approximation obtained from the TSFP \eqref{psifull}--\eqref{ufull} (or equivalently \eqref{ufull2}). Under the assumption (A), there exist
$h_0 > 0$ and $\tau_0 > 0$ sufficiently small and independent of $\varepsilon$ such that,
for any $0 < \varepsilon \leq 1$, when $0 < h \leq h_0$ and $0< \tau \leq \tau_0\varepsilon^{(\beta-p)/2}$, we have the error estimates for $s\in (1/2, m]$
\begin{equation}
\|u(\cdot, t_n) - I_N(u^n)\|_s+\|\partial_t u(\cdot, t_n) - I_N( v^n)\|_{s-1} \lesssim h^{1+m-s} + \varepsilon^{p-\beta}\tau^2,\quad 0 \leq n \leq \frac{T_0/\varepsilon^\beta}{\tau}.
\label{eq:error1}
\end{equation}
Furthermore, there exists a constant $M>0$ depending on $T_0$, $\|u_0\|_{m+1}$, $\|u_1\|_m$, $\|u\|_{L^\infty([0, T_\varepsilon]; H_{\rm per}^{m})}$ and  $\|\partial_t u\|_{L^\infty([0, T_\varepsilon]; H_{\rm per}^{m-1})}$ such that the numerical solution satisfies
\be\label{reg}
\|I_N(u^n)\|_{m+1}+\|I_N( v^n)\|_m\le M,\quad 0 \leq n \leq \frac{T_0/\varepsilon^\beta}{\tau}.
\ee
\end{theorem}

\begin{remark}
It follows from the $\varepsilon$-dependent error estimate that large time step at $\tau\sim \varepsilon^{-(p-\beta)/2}$ when $0\le \beta<p$ is allowed to simulate the long-time dynamics of the NKGE up to time $T_0/\varepsilon^\beta$. Particularly, the error bound is uniformly at the second order for the large time step $\tau=O(\varepsilon^{-(p-\beta)/2})$ in the parameter regime $0\le \beta <p$. While for $\beta=p$, $\tau$ has to be taken as $O(1)$.
\end{remark}

\begin{remark}
The results in Theorem \ref{thm:eb1} are still valid in high dimensions, i.e., $d > 1$, if $m >\frac{d}{2}$.
\end{remark}

\subsection{Preliminary estimates}
In this subsection, we prepare some results for proving the main theorem.
Denote
\[F_t: \ \phi \mapsto e^{-it\langle\nabla\rangle}  F\big(e^{it\langle\nabla\rangle}\phi\big),\quad t \in \mathbb{R},\]
where $F$ is defined by \eqref{Fd}, then we have the following proposition on the properties of $F_t$.
\begin{proposition}
(i) Let $s>1/2$, then for any $t \in \mathbb{R}$, the function $F_t: H^{s}(\Omega) \to H^{s+1}(\Omega)$ is $C^{\infty}$ and satisfies
\be\label{Fb}
\begin{split}
&\left\|F_t(\phi)\right\|_{s+1}\le C \|\phi\|_{s}^{p+1},\quad
\left\|F'_t(\phi)(\eta)\right\|_{s+1}\le C \|\phi\|_{s}^{p}\,\|\eta\|_{s},\\
&\left\|F''_t(\phi)(\eta, \zeta)\right\|_{s+1} \leq C\|\phi\|^{p-1}_{s}\,\|\eta\|_{s}\,\|\zeta\|_{s}.
\end{split}
\ee
(ii) If $s\ge 1$, then the derivatives with respect to $t$ satisfy
\be\label{Fdb}
\left\|\partial_t F_t(\phi)\right\|_{s}\le C\|\phi\|_{s}^{p+1},\,\,\,\,
\left\|\partial^2_t F_t(\phi)\right\|_{s} \leq C\|\phi\|_{s+1}^{p+1},\,\,\,\,
\left\|\partial_t F'_t(\phi)(\eta)\right\|_{s} \leq C\|\phi\|_{s}^{p}\|\eta\|_{s}.
\ee
(iii) Assume $s>1/2$, $\phi, \eta\in B_R^s:=\{v\in H^s(\Omega), \|v\|_s\le R\}$, then there exists a constant $L>0$ depending on $R$ such that for all $t \in \mathbb{R}$ and $\sigma\in[0, s]$, the Lipschitz estimate is valid:
\be\label{Fl}
\|G(\phi) -G(\eta)\|_{\sigma} \leq L \|\phi-\eta\|_{\sigma},\quad
\|F_t(\phi) - F_t(\eta)\|_{\sigma+1} \leq L \|\phi-\eta\|_{\sigma}.
\ee
\end{proposition}
\emph{Proof.}
Firstly, we recall the inequality which was established in \cite{CMTZ}:
\be\label{tame}
\|vw\|_\sigma\le C\|v\|_\sigma\,\|w\|_s, \quad v\in H^\sigma(\Omega),\quad w\in H^s(\Omega),
\ee
for $s>1/2$ and $\sigma\in [0, s]$. Hence for $\phi\in H^s(\Omega)$, one has
\begin{align*}
\left\|F_t(\phi)\right\|_{s+1}&=\left\|F\left(e^{it\langle\nabla\rangle}
\phi\right)\right\|_{s+1}
=\left\|f\left(\frac{1}{2}\left(e^{it\langle\nabla\rangle}\phi + e^{-it\langle\nabla\rangle}\overline{\phi}\right)\right)\right\|_{s}\\
&\le C\left\|e^{it\langle\nabla\rangle}\phi + e^{-it\langle\nabla\rangle}\overline{\phi}\right\|_s^{p+1}\le C\|\phi\|_s^{p+1}.
\end{align*}
Noticing that $f(v)=v^{p+1}$, a direct calculation gives
\be\label{Fp}
F'(\phi)(\eta)=\frac{(p+1)i}{2^{p+1}}\langle\nabla\rangle^{-1}\left((\phi + \overline{\phi})^{p}(\eta + \overline{\eta})\right),
\ee
which implies that
\be \label{Fpn}
\|F'(\phi)(\eta)\|_{s+1}=\frac{p+1}{2^{p+1}}\left\|(\phi + \overline{\phi})^{p}(\eta + \overline{\eta})\right\|_{s}\le C\|\phi\|_s^{p}\|\eta\|_s.
\ee
Note that
\[F_t'(\phi)(\eta)=e^{-it\langle\nabla\rangle}F'
\left(e^{it\langle\nabla\rangle}\phi\right)
\left(e^{it\langle\nabla\rangle}\eta\right),\]
and this immediately yields the second inequality in \eqref{Fb}. The second derivative of $F$
takes the form
\[F''(\phi)(\eta, \zeta)=\frac{p(p+1)i}{2^{p+1}}\langle\nabla\rangle^{-1}\left((\phi + \overline{\phi})^{p-1}(\eta + \overline{\eta})(\zeta + \overline{\zeta})\right),\]
which leads to that
\[\left\|F''(\phi)(\eta, \zeta)\right\|_{s+1}=\frac{p(p+1)}{2^{p+1}}\left\|(\phi + \overline{\phi})^{p-1}(\eta + \overline{\eta})(\zeta + \overline{\zeta})\right\|_s\le C\|\phi\|_s^{p-1}\|\eta\|_s\|\zeta\|_s.\]
Thus the last inequality in \eqref{Fb} can be obtained by recalling
\[F_t''(\phi)(\eta, \zeta)=e^{-it\langle\nabla\rangle}F''\left(e^{it\langle\nabla\rangle}\phi\right)
\left(e^{it\langle\nabla\rangle}\eta, e^{it\langle\nabla\rangle}\zeta\right).\]

The first derivative of $F_t$ with respect to $t$ reads as
\[\partial_t F_t(\phi)=-i\langle\nabla\rangle F_t(\phi)+e^{-it\langle\nabla\rangle}  F'(\mu)(i \langle\nabla\rangle\mu),\quad \mu=e^{it\langle\nabla\rangle}\phi.\]
Applying \eqref{Fb}, \eqref{tame} and \eqref{Fp}, we obtain
\begin{align*}
\|\partial_t F_t(\phi)\|_s&\le \|F_t(\phi)\|_{s+1}+\|F'(\mu)(i\langle\nabla\rangle\mu)\|_s\\
&\le C\|\phi\|_s^{p+1}+ C\|(\mu+\overline{\mu})^{p}(\langle\nabla\rangle\mu-\langle\nabla\rangle\overline{\mu})\|_{s-1}\\
&\le C\|\phi\|_s^{p+1}+ C\|(\mu+\overline{\mu})^{p}\|_s\|\langle\nabla\rangle(\mu-\overline{\mu})\|_{s-1}\\
&\le C\|\phi\|_s^{p+1}+ C\|\mu+\overline{\mu}\|^{p}_s\|\mu-\overline{\mu}\|_{s}\\
&\le C\|\phi\|_s^{p+1}.
\end{align*}
Further computations give that
\begin{align*}
\partial_t^2 F_t(\phi)&=-\langle\nabla\rangle^2F_t(\phi)-2i\langle\nabla\rangle
e^{-it\langle\nabla\rangle}F'(\mu)(i\langle\nabla\rangle\mu)+e^{-it\langle\nabla\rangle}
F'(\mu)(-\langle\nabla\rangle^2\mu)\\
&\quad+e^{-it\langle\nabla\rangle}F''(\mu)(i\langle\nabla\rangle\mu,
i\langle\nabla\rangle\mu),\end{align*}
which leads to
\begin{align*}
\|\partial_t^2 F_t(\phi)\|_s&\le
\|F_t(\phi)\|_{s+2}+2\|F'(\mu)(i\langle\nabla\rangle\mu)\|_{s+1}+
\|F'(\mu)(-\langle\nabla\rangle^2\mu)\|_s\\
&\quad+\|F''(\mu)(i\langle\nabla\rangle\mu,
i\langle\nabla\rangle\mu)\|_s\\
&\hspace{-3mm}\le C\|\phi\|_{s+1}^{p+1}+C\|(\mu+\overline{\mu})^{p}
\langle\nabla\rangle^2(\mu+\overline{\mu})\|_{s-1}+
C\|(\mu + \overline{\mu})^{p-1}(\langle\nabla\rangle(\mu-\overline{\mu}))^{2}\|_{s-1}\\
&\hspace{-3mm}\le C\|\phi\|_{s+1}^{p+1}+C \|\mu+\overline{\mu}\|_s^{p}\|\mu+\overline{\mu}\|_{s+1}+
C\|\mu+\overline{\mu}\|_s^{p-1}\|\mu-\overline{\mu}\|_{s+1}^{2}\\
&\hspace{-3mm}\le C\|\phi\|_{s+1}^{p+1}.
\end{align*}
For the last inequality of \eqref{Fdb}, note that
\[\partial_t F_t'(\phi)(\eta)=-i\langle\nabla\rangle F_t'(\phi)(\eta)+
e^{-it\langle\nabla\rangle}F''(\mu) (\nu, i\langle\nabla\rangle\mu)+
e^{-it\langle\nabla\rangle}F'(\mu)(i\langle\nabla\rangle\nu),\]
where $\nu=e^{it\langle\nabla\rangle}\eta$.
Thus we get
\begin{align*}
\|\partial_t F_t'(\phi)(\eta)\|_s&\le \|F_t'(\phi)(\eta)\|_{s+1}+
\|F''(\mu) (\nu, i\langle\nabla\rangle\mu)\|_s+\|F'(\mu)(i\langle\nabla\rangle\nu)\|_s\\
&\le C\|\phi\|_s^{p}\|\eta\|_s+ C\|(\mu + \overline{\mu})^{p-1}(\nu + \overline{\nu})\langle\nabla\rangle (\mu-\overline{\mu})\|_{s-1}+
C\|\mu+\overline{\mu}\|^{p}_s\|\nu-\overline{\nu}\|_{s}\\
&\le C\|\phi\|_s^{p}\|\eta\|_s+ C\|\mu + \overline{\mu}\|_s^{p-1}\|\nu + \overline{\nu}\|_s\|\langle\nabla\rangle \mu-\langle\nabla\rangle\overline{\mu}\|_{s-1}\\
&\le C\|\phi\|_s^{p}\|\eta\|_s,
\end{align*}
which completes the proof for \eqref{Fdb}.

For the Lipschitz estimate \eqref{Fl}, a straightforward calculation shows that
\begin{align*}
\|G(\phi)-G(\eta)\|_{\sigma}&=\Big\|f\big(\frac{1}{2}(\phi+\overline{\phi})\big)-
f\big(\frac{1}{2}(\eta+\overline{\eta})\big)\Big\|_\sigma\\
&=\frac{1}{2^{p+1}}\left\|\left[\sum_{q=0}^{p}\binom{p+1}{q}
(\phi+\overline{\phi}-\eta-\overline{\eta})^{p-q}(\eta+\overline{\eta})^q
\right](\phi-\eta+\overline{\phi}-\overline{\eta})\right\|_\sigma\\
&\le \frac{1}{2^{p+1}}\sum_{q=0}^{p}\binom{p+1}{q}
\left\|\phi+\overline{\phi}-\eta-\overline{\eta}\right\|_s^{p-q}\left\|
\eta+\overline{\eta}\right\|_s^q\left\|\phi-\eta+\overline{\phi}-\overline{\eta}\right\|_\sigma\\
&\le CR^{p}\left\|\phi-\eta\right\|_\sigma.
\end{align*}
Noticing that
\begin{align*}
\|F_t(\phi)-F_t(\eta)\|_{\sigma+1}&=\big\|F(e^{it \langle\nabla\rangle}\phi)-F(e^{it\langle\nabla\rangle}\eta)\big\|_{\sigma+1}=\big\|G(e^{it \langle\nabla\rangle}\phi)-G(e^{it\langle\nabla\rangle}\eta)\big\|_{\sigma}\\
&\le CR^{p}\left\|\phi-\eta\right\|_\sigma,
\end{align*}
the proof is completed.
\hfill $\square$ \bigskip

Concerning on the flow $\mathcal{S}_\tau$ in \eqref{Strang}, we have the stability estimate as follows.
\begin{lemma}\label{stab}
Assume $\phi_0, \eta_0\in B_R^s$ with $s>1/2$, then for any $\tau>0$, we have
\[
\|\mathcal{S}_{\tau} (\phi_0) -\mathcal{S}_{\tau} (\eta_0) \|_s \leq (1+L\varepsilon^{p}  \tau)\|\phi_0-\eta_0\|_s,
\]
where $L$ depends on $R$.
\end{lemma}
\emph{Proof.}
Since the operator $e^{i\tau\langle\nabla\rangle}$ is an isometry, we only need to consider the operator associated with the nonlinear subproblem.
By the definition and the Lipschitz estimate \eqref{Fl}, we have
\begin{align*}
\|\varphi^{\tau}_V(\phi_0)-\varphi^{\tau}_V(\eta_0)\|_s &\leq \|\phi_0- \eta_0\|_s+\varepsilon^{p} \tau\|F(\phi_0)-F(\eta_0)\|_s\\
&\leq \|\phi_0-\eta_0\|_s+L \varepsilon^{p} \tau\|\phi_0-\eta_0\|_{s}\\
&\le (1+L\varepsilon^{p}  \tau)\|\phi_0-\eta_0\|_s,
\end{align*}
which completes the proof.
\hfill $\square$ \bigskip

\begin{lemma}\label{localeror}
Denote the exact solution of \eqref{eq:NLS} with initial data $\psi_0$ as
$\psi(t)=\mathcal{S}_{e,t}(\psi_0)$. Assume $\psi(t)\in H^{s+1} (s\ge 1)$, then for $0<\varepsilon\le 1$ and $0<\tau\le 1/\varepsilon^{p}$, the local error of the Strang splitting \eqref{Strang} is bounded by
\[\|\mathcal{S}_\tau(\psi(t_n))-\mathcal{S}_{e,\tau}(\psi(t_n))\|_s\le M_0\varepsilon^{p}\tau^3,\]
where $M_0$ depends on $\|\psi\|_{L^\infty([0,T_\varepsilon]; H^{s+1})}$.
\end{lemma}
\emph{Proof.}
For simplicity of notation, we denote $\psi_n=\psi(t_n)$.
An application of the Duhamel's principle leads to the following representation of the exact solution
\begin{equation}
\psi(t_n+t)=e^{it\langle\nabla\rangle}\psi_n+\varepsilon^{p} e^{it\langle\nabla\rangle}
\int^{t}_0 e^{-i\theta\langle\nabla\rangle}F\left( \psi(t_n+\theta)\right) d\theta.
\end{equation}
Introducing $\eta_n(t):=e^{-i(t_n+t)\langle\nabla\rangle}\psi(t_n+t)$, we have
\begin{equation}
\eta_n(t)=\eta_n(0)+\varepsilon^{p}  \int^{t}_0 F_{t_n+\theta}(\eta_n(\theta)) d\theta.
\end{equation}
Applying the Taylor expansion
\[F_t(z_1+z_2)=F_t(z_1)+F_t'(z_1)(z_2)+\int_0^1 (1-\theta)F_t''(z_1+\theta z_2)(z_2^2)d\theta,\]
we yield
\begin{align*}
\eta_n(\tau)&=\eta_n(0)+\varepsilon^{p}  \int^{\tau}_0 F_{t_n+\theta}\Big(\eta_n(0)+\varepsilon^{p}\int^\theta_0 F_{t_n+\theta_1}\left(\eta_n(\theta_1) \right) d\theta_1 \Big) d\theta\\
&=\eta_n(0)+\varepsilon^{p}\int^{\tau}_0 F_{t_n+\theta}(\eta_n(0)) d\theta+\varepsilon^{2p}\int^{\tau}_0 \int^\theta_0 F'_{t_n+\theta}(\eta_n(0))F_{t_n+\theta_1}(\eta_n(\theta_1))d\theta_1 d\theta \\
&\,\,\,+\varepsilon^{3p} \int^{\tau}_0\int^1_0(1-\zeta)F''_{t_n+\theta}\left((1-\zeta)\eta_n(0) + \zeta \eta_n(\theta)\right)\Big(\int^\theta_0F_{t_n+\theta_1} (\eta_n(\theta_1))d\theta_1\Big)^2 d\zeta d\theta\\
&=\eta_n(0)+\varepsilon^{p} \int^{\tau}_0 F_{t_n+\theta}(\eta_n(0)) d\theta+\varepsilon^{2p} \int^{\tau}_0 \int^\theta_0 F'_{t_n+\theta}(\eta_n(0))F_{t_n+\theta_1}(\eta_n(0))d\theta_1 d\theta\\
&\,\,\,+\varepsilon^{3p} \int^{\tau}_0\int^1_0(1-\zeta)F''_{t_n+\theta}\left((1-\zeta)\eta_n(0)+ \zeta \eta_n(\theta)\right)\Big(\int^\theta_0F_{t_n+\theta_1} \left(\eta_n(\theta_1)\right)d\theta_1\Big)^2 d\zeta d\theta\\
&\,\,\,+\varepsilon^{3p} \int^{\tau}_0\int^\theta_0\int^1_0 F'_{t_n+\theta}(\eta_n(0))F'_{t_n+\theta_1} \left((1-\zeta)\eta_n(0)+\zeta\eta_n(\theta_1)\right)\\
&\qquad\qquad\qquad\qquad\qquad\qquad\qquad\qquad\qquad\qquad
\Big(\int^{\theta_1}_0F_{t_n+\theta_2}
\left(\eta_n(\theta_2)\right)d\theta_2\Big)d\zeta d\theta_1 d\theta.
\end{align*}
Twisting the variable back, we obtain
\begin{align}
\mathcal{S}_{e,\tau}(\psi_n)&=e^{i(t_n+\tau)\langle\nabla\rangle}\eta_n(\tau)\nn\\
&=e^{i\tau\langle\nabla\rangle}\psi_n+ \varepsilon^{p} e^{i\tau\langle\nabla\rangle}\int^{\tau}_0 F_\theta\left( \psi_n\right)d\theta  +\varepsilon^{3p}  e^{i\tau\langle\nabla\rangle} E_3\nn\\
&\quad+ \varepsilon^{2p} e^{i\tau\langle\nabla\rangle}\int^{\tau}_0 \int^\theta_0 F'_\theta\left(\psi_n\right) F_{\theta_1} \left(\psi_n\right) d\theta_1 d\theta,\label{exact}
\end{align}
where $E_{3}= E_{3, 1}+ E_{3, 2}$ with
\begin{align*}
&E_{3, 1}= \int^{\tau}_0\int^1_0(1-\zeta)F''_\theta\left((1-\zeta)\psi_n+\zeta e^{-i\theta\langle\nabla\rangle} \psi(t_n+\theta)\right)\\
&\qquad\qquad\qquad\qquad\qquad\qquad\qquad\qquad\Big(\int^\theta_0 F_{\theta_1} \Big(e^{-i\theta_1\langle\nabla\rangle}\psi(t_n+\theta_1)\Big)d\theta_1\Big)^2 d\zeta d\theta,\\
&E_{3, 2}=\int^{\tau}_0 \int^\theta_0\int^1_0 F'_\theta(\psi_n) F'_{\theta_1} \big((1-\zeta)\psi_n+\zeta e^{-i\theta_1\langle\nabla\rangle}\psi(t_n+\theta_1)\big)\\
&\qquad\qquad\qquad\qquad\qquad\qquad\qquad\qquad\Big(\int^{\theta_1}_0 F_{\theta_2}\big(e^{-i\theta_2\langle\nabla\rangle}\psi(t_n+\theta_2)\big)d\theta_2\Big)d\zeta d\theta_1 d\theta.
\end{align*}

On the other hand, noticing \eqref{nlex1}, for the Strang splitting we get
\[
\mathcal{S}_{\tau}(\psi_n)= e^{i\tau\langle\nabla\rangle/2}\left[e^{i\tau\langle\nabla\rangle/2}\psi_n+
\varepsilon^{p} \tau F\big(e^{i\tau\langle\nabla\rangle/2}\psi_n\big)\right]=e^{i\tau\langle\nabla\rangle}
\left(\psi_n+\varepsilon^{p}\tau F_{\tau/2}(\psi_n)\right).
\]
Then the local truncation error can be written as
\be\label{locl}
\mathcal{S}_{\tau}(\psi_n) -\mathcal{S}_{e, \tau}(\psi_n)
=\varepsilon^{p}  e^{i\tau\langle\nabla\rangle}r_1-\varepsilon^{2p} e^{i\tau\langle\nabla\rangle} r_2-\varepsilon^{3p} e^{i\tau\langle\nabla\rangle}  E_3,
\ee
where
\[r_1=\tau F_{\tau/2}(\psi_n)-\int^{\tau}_0 F_\theta(\psi_n)d\theta,\quad
r_2=\int^{\tau}_0 \int^\theta_0 F'_\theta\left(\psi_n\right) F_{\theta_1} \left(\psi_n\right) d\theta_1 d\theta.\]
Next we estimate each term individually.
Express the quadrature error in the second-order Peano form,
\[r_1=-\tau^3\int^1_0 \kappa_2(\theta) \partial^2_{\omega} F_{\omega}(\psi_n) |_{\omega=\theta \tau} d\theta,\quad {\kappa_2(\theta)=\frac{1}{2}\min\{\theta^2, (1-\theta)^2\}}.\]
Applying \eqref{Fdb}, we obtain
\be\label{r1e}
\|r_1\|_s\le C\tau^3\, \|\psi_n\|_{s+1}^{p+1}\int^1_0 \kappa_2(\theta)d\theta\lesssim\tau^3.
\ee
Inserting the identities
\[
F_{\theta_1}(\psi_n) = F_{\tau/2}(\psi_n) + \int^{\theta_1}_{\tau/2} \partial_{\omega} F_{\omega}(\psi_n) d\omega, \quad F'_\theta(\psi_n)=F'_{\tau/2}(\psi_n)+\int^{\theta}_{\tau/2} \partial_{\omega}F'_{\omega}(\psi_n)d\omega\]
into the double integral term, we get
\begin{align*}
r_2 &=\frac{1}{2}\tau^2F'_{\tau/2}(\psi_n)F_{\tau/2}(\psi_n)+\int^{\tau}_0 \int^{\theta}_0 F'_{\tau/2}(\psi_n) \int^{\theta_1}_{\tau/2} \partial_{\omega} F_{\omega} (\psi_n)d\omega d \theta_1 d\theta\\
&\quad+\int^{\tau}_0 \theta\int^{\theta}_{\tau/2} \,\partial_{\omega}F'_{\omega}(\psi_n) F_{\tau/2}(\psi_n) d\omega d\theta+\int^{\tau}_0 \int^{\theta}_0 \int^{\theta}_{\tau/2} \int^{\theta_1}_{\tau/2} \partial_{\omega}F'_{\omega}(\psi_n) \partial_{\omega_1} F_{\omega_1}(\psi_n) d\omega_1 d\omega  d\theta_1 d\theta.
\end{align*}
By definition, we have
\[F'_{\tau/2}(\psi_n)F_{\tau/2}(\psi_n)=e^{-i\frac{\tau}{2}\langle \nabla \rangle}F'(e^{i\frac{\tau}{2}\langle \nabla \rangle}\psi_n)\big(F(e^{i\frac{\tau}{2}\langle \nabla \rangle}\psi_n)\big)=0,\]
by recalling \eqref{Fp} and the fact that $F(\cdot)$ is purely imaginary. Applying \eqref{Fb} and \eqref{Fdb}, we obtain
\be\label{r2e}
\begin{split}
\|r_2\|_s &\le C\tau^3 \|\psi_n\|_s^{p} \sup\limits_{0\le \omega\le\tau}\left\|\partial_{\omega} F_{\omega} (\psi_n)\right\|_s+C\tau^3\|\psi_n\|_s^{p}
\left\|F_{\tau/2}(\psi_n)\right\|_s\\
&\quad+C\tau^4\|\psi_n\|_s^{p}\sup\limits_{0\le \omega\le\tau}\left\|\partial_{\omega} F_{\omega} (\psi_n)\right\|_s\\
&\le C(\tau^3+\tau^4)\|\psi_n\|_s^{2p+1}\lesssim (\tau^3+\tau^4).
\end{split}
\ee
Using \eqref{Fb}, we derive
\begin{align}
\|E_3\|_s&\le \|E_{3,1}\|_s+\|E_{3, 2}\|_s\nn\\
&\le C\tau^3 \sup\limits_{0\le \theta\le \tau}\|\psi(t_n+\theta)\|_s^{p-1} \sup\limits_{0\le \theta\le\tau}\left\|F_\theta \big(e^{-i\theta\langle\nabla\rangle}\psi(t_n+\theta)\big)\right\|_s^2\nn\\
&\quad+C\tau^3\|\psi_n\|_s^{p} \sup\limits_{0\le \theta\le \tau}\|\psi(t_n+\theta)\|_s^{p}\sup\limits_{0\le \theta\le\tau}\left\|F_\theta \big(e^{-i\theta\langle\nabla\rangle}\psi(t_n+\theta)\big)\right\|_s\nn\\
&\le C\tau^3\sup\limits_{0\le \theta\le \tau}\|\psi(t_n+\theta)\|_s^{3p+1}\lesssim \tau^3.\label{e3e}
\end{align}
Combining \eqref{locl}-\eqref{e3e}, we arrive at the conclusion and the proof is complete.
\hfill $\square$ \bigskip

\subsection{Proof of Theorem \ref{thm:eb1} }
Similar to the proof of the TSFP method for the Dirac equation \cite{BCJY}, the proof will be divided into two parts: (I) to prove the convergence of the semi-discretization, and (II) to complete the error analysis by comparing the semi-discretization \eqref{Strang} and the full-discretization \eqref{psifull}.

\smallskip
\noindent
{\bf{Part I}} (Convergence of the semi-discretization) Firstly, we observe that the assumption (A) is equivalent to the regularity of $\psi(x, t)$ as
\[\psi\in L^\infty\left([0, T_\varepsilon]; H_{\rm per}^{m+1}\right),\quad
\|\psi\|_{L^{\infty}\left([0, T_{\varepsilon}]; H^{m+1}_{\rm per}\right)} \lesssim 1.\]
Now, we give a global error on the Strang splitting \eqref{Strang}: there exists $\tau_0>0$ independent of $\varepsilon$ such that when $0 < \tau\le \tau_0\varepsilon^{(\beta-p)/2}$, the error of the Strang splitting satisfies
\be
\label{semi-eror}
\|\psi^{[n]}-\psi(\cdot,t_n)\|_m\le M_1\varepsilon^{p-\beta}\tau^2,\quad \|\psi^{[n]}\|_m\le R+1,\quad
0\le n\le \frac{T_0/\varepsilon^\beta}{\tau},
\ee
where $R:=\|\psi\|_{L^\infty([0, T_\varepsilon]; H_{\rm per}^{m})}$ and $M_1$ depends on $T_0$, $R$ and $\|\psi\|_{L^\infty([0, T_\varepsilon]; H_{\rm per}^{m+1})}$. Furthermore, for the regularity of $\psi^{[n]}$, we have $\psi^{[n]}\in H_{\rm per}^{m+1}$ when $\tau\le \tau_0\varepsilon^{(\beta-p)/2}$ with
\be
\label{semi-re}
\|\psi^{[n]}\|_{m+1}\le M_2, \quad 0\le n\le \frac{T_0/\varepsilon^\beta}{\tau},
\ee
where $M_2$ depends on $T_0$, $R$ and $\|\psi_0\|_{m+1}$.

We apply a standard induction argument for proving \eqref{semi-eror}. Firstly, it is obvious for $n=0$ since $\psi^{[0]}=\psi_0\in B^{m}_R$. Assume $\psi^{[k]}\in B^{m}_{R+1}$ for $0\leq k\leq n<\frac{T_0/\varepsilon^\beta}{\tau}$. Denote $e^{[k]}=\psi^{[k]}-\psi(\cdot, t_k)$.
By definition,
\[e^{[k+1]}=\mathcal{S}_\tau(\psi^{[k]})-\mathcal{S}_\tau(\psi(t_k))+
\mathcal{S}_\tau(\psi(t_k))-\mathcal{S}_{e,\tau}(\psi(t_k)).\]
Using Lemmas \ref{stab} and \ref{localeror}, we get when $\tau\le 1/\varepsilon^{p}$,
\[\big\|e^{[k+1]}\big\|_{m}-\|e^{[k]}\|_{m} \leq L\varepsilon^{p} \tau \big\|e^{[k]}\big\|_{m}+M_0\varepsilon^{p}\tau^3,\]
where $L$ and $M_0$ depend on $R$  and $\big\|\psi\big\|_{L^\infty([0, T_\varepsilon]; H_{\rm per}^{m+1})}$, respectively,  as claimed in Lemmas \ref{stab} and \ref{localeror}.
Summing the above inequality for $k=0, \ldots,n$, one gets
\begin{align*}
\big\|e^{[n+1]}\big\|_{m}&\leq \big\|e^{[0]}\big\|_{m}+L\varepsilon^{p} \tau \sum\limits_{k=0}^n \big\|e^{[k]}\big\|_{m}+M_0\varepsilon^{p}\tau^3(n+1)\\
&\le M_0T_0\varepsilon^{p-\beta}\tau^2+L\varepsilon^{p}  \tau \sum\limits_{k=0}^n \big\|e^{[k]}\big\|_{m}.
\end{align*}
 Applying the Gronwall's inequality, we derive
\[\big\|e^{[n+1]}\big\|_{m}\le M_0T_0 e^{LT_0} \varepsilon^{p-\beta}\tau^2,\quad 0\le n< \frac{T_0/\varepsilon^\beta}{\tau}.\]
Then the triangle inequality yields that
\[\big\|\psi^{[n+1]}\big\|_{m}\leq\big\|\psi(\cdot, t_{n+1})\big\|_{m}+1, \quad 0\leq n< \frac{T_0/\varepsilon^\beta}{\tau},\]
when $0<\tau\leq 1/\varepsilon^{p}$ and $\tau\le (M_0T_0)^{-1/2}e^{-LT_0/2}\varepsilon^{(\beta-p)/2}$. Set $\tau_0=\min\{1, (M_0T_0)^{-1/2}e^{-LT_0/2}\}$, then the induction \eqref{semi-eror} holds when $\tau\le \tau_0\varepsilon^{(\beta-p)/2}$ and $\varepsilon\in (0, 1]$.
For the last inequality \eqref{semi-re}, recalling \eqref{nlex1} and \eqref{Fb}, we have
\begin{align*}
\big\|\psi^{[n+1]}\big\|_{m+1}&=\big\|\varphi_V^\tau(e^{i\tau/2\langle \nabla \rangle }\psi^{[n]})\big\|_{m+1}\\
&\le \big\|e^{i\tau/2\langle \nabla \rangle }\psi^{[n]}\big\|_{m+1}+\varepsilon^{p}\tau\left\|F\big(e^{i\tau/2\langle \nabla \rangle }\psi^{[n]}\big)\right\|_{m+1}\\
&\le \big\|\psi^{[n]}\big\|_{m+1}+C\varepsilon^{p}\tau \big\|\psi^{[n]}\big\|^{p+1}_{m}\\
&\le \big\|\psi^{[n]}\big\|_{m+1}+C\varepsilon^{p}\tau (R+1)^{p+1}\\
&\le \big\|\psi^{[0]}\big\|_{m+1}+C(n+1)\varepsilon^{p}\tau (R+1)^{p+1}\\
&\le \big\|\psi_0\big\|_{m+1}+CT_0 (R+1)^{p+1},
\end{align*}
and \eqref{semi-re} is established.

\smallskip

\noindent
{\bf{Part II}} (Convergence of the full-discretization) For $0 \leq n \leq \frac{T_0/\varepsilon^\beta}{\tau}$, we rewrite the error as
\be\label{sp}
\psi(\cdot, t_n)-I_N(\psi^n)=\psi(\cdot, t_n)-\psi^{[n]}+\psi^{[n]}-P_N(\psi^{[n]})+ P_N(\psi^{[n]})-I_N(\psi^n).
\ee
For $0\le s\le m$, the regularity result \eqref{semi-re} implies that
\be\label{psip}
\|\psi^{[n]}-P_N(\psi^{[n]})\|_s\leq C M_2 h^{1+m-s},
\ee
and by \eqref{semi-eror},
\be\label{psi23}
\|\psi(\cdot, t_n) - \psi^{[n]}\|_s \le \|\psi(\cdot, t_n) - \psi^{[n]}\|_m\leq M_1\varepsilon^{p-\beta}\tau^2.
\ee
Thus, it remains to establish the error bound for the error
\[
e^n:=P_N(\psi^{[n]})-I_N(\psi^n),\quad 0 \leq n \leq \frac{T_0/\varepsilon^\beta}{\tau}.\]
Now, we'll use an induction to show that when $h$ is sufficiently small, we have
\be\label{eq:psib}
\|e^n\|_l\leq  M_3 h^{1+m-l},\,\,\,\,
l\in (1/2, m+1];\,\,\,\, \|I_N(\psi^n)\|_m\leq C(1+R)+1,
\ee
where $M_3$ depends on $T_0$, $R$ and $\|\psi_0\|_{m+1}$.

For $n=0$, \eqref{eq:psib} is obvious by using the projection and interpolation errors \cite{ST}:
\begin{align*}
&\|e^0\|_l=\|P_N(\psi_0)-I_N(\psi_0)\|_l\le C_1h^{1+m-l}\|\psi_0\|_{m+1},\\
&\|I_N(\psi^0)\|_m\le \|P_N(\psi_0)\|_m+\|e^0\|_m\le C\|\psi^0\|_m+C_1h\|\psi_0\|_{m+1}\le C(1+R)+1,
\end{align*}
when $h\le \frac{1}{C_1\|\psi_0\|_{m+1}}$.
For $n\ge 1$, assume \eqref{eq:psib} holds for $0\leq k\leq n<\frac{T_0/\varepsilon^\beta}{\tau}$. We rewrite \eqref{psifull} as
\begin{align*}
&\psi^{(n,1)}=e^{i\tau\langle\nabla\rangle/2}I_N(\psi^n),\quad
\psi^{(n,2)}=\psi^{(n,1)}+i\varepsilon^{p} \tau \langle \nabla \rangle^{-1}I_N(G(\psi^{(n,1)})),\\
&I_N(\psi^{n+1})=e^{i\tau\langle\nabla\rangle/2}I_N(\psi^{(n,2)}).
\end{align*}
Hence we get $\psi^{(n,1)}, \psi^{(n,2)}\in Y_N$.
Similarly, \eqref{Strang} can be expressed as
\[\psi^{\langle n,1\rangle}=e^{i\tau\langle\nabla\rangle/2}\psi^{[n]},\,\,\,
\psi^{\langle n, 2\rangle}=\psi^{\langle n, 1\rangle}+i\varepsilon^{p} \tau \langle \nabla \rangle^{-1} G(\psi^{\langle n, 1\rangle}),\,\,\,
\psi^{[n+1]}=e^{i\tau\langle\nabla\rangle/2}\psi^{\langle n, 2\rangle},\]
which implies that
\begin{align*}
&P_N(\psi^{\langle n,1\rangle})=e^{i\tau\langle \nabla \rangle/2}P_N(\psi^{[n]}),\\
&P_N(\psi^{\langle n,2\rangle})=P_N(\psi^{\langle n,1\rangle})+i\varepsilon^{p} \tau \langle \nabla \rangle^{-1} P_N(G(\psi^{\langle n,1\rangle})),\\
&P_N(\psi^{[n+1]})=e^{i\tau\langle\nabla\rangle/2}P_N(\psi^{\langle n, 2\rangle}).
\end{align*}
Thus by definition, we get
\begin{align*}
&\hspace{-4mm}\|e^{n+1}\|_l=\big\|P_N(\psi^{[n+1]})-I_N(\psi^{n+1})\big\|_l=\big\|
P_N(\psi^{\langle n, 2\rangle})-
I_N(\psi^{(n,2)})\big\|_l\\
&\le \big\|P_N(\psi^{\langle n, 1\rangle})-I_N(\psi^{(n,1)})\big\|_l+\varepsilon^{p}\tau \big\|P_N(G(\psi^{\langle n,1\rangle}))-I_N(G(\psi^{(n,1)}))\big\|_{l-1}\\
&\le \|e^n\|_l+\varepsilon^{p}\tau \big\|P_N(G(\psi^{\langle n, 1\rangle}))-I_N(G(\psi^{\langle n, 1\rangle}))\big\|_l\\
&\quad+\varepsilon^{p}\tau \big\|I_N(G(\psi^{\langle n, 1\rangle}))-I_N(G(\psi^{(n,1)}))\big\|_{\min\{l,m\}}\\
&\le \|e^n\|_l+C\varepsilon^{p}\tau h^{1+m-l} \big\|G(\psi^{\langle n, 1\rangle})\big\|_{m+1} +C\varepsilon^{p}\tau\big\|G(\psi^{\langle n, 1\rangle})-G(\psi^{(n,1)})\big\|_{\min\{l,m\}}\\
&\le \|e^n\|_l+C\varepsilon^{p}\tau h^{1+m-l}\|\psi^{\langle n, 1\rangle}\|^{p+1}_{m+1}+CL\varepsilon^{p}\tau\|\psi^{\langle n, 1\rangle}-\psi^{(n,1)}\|_l\\
&\le (1+CL\varepsilon^{p}\tau)\|e^n\|_l+CM_2^{p+1}
\varepsilon^{p}\tau h^{1+m-l}+CL\varepsilon^{p}\tau \big\|P_N(\psi^{[n]})-\psi^{[n]}\big\|_l\\
&\le (1+CL\varepsilon^{p}\tau)\|e^n\|_l+CM_2(L+ M_2^{p})
\varepsilon^{p}\tau h^{1+m-l},
\end{align*}
where we have used the fact that $\psi^{[n]}, \psi^{\langle n,1\rangle}, G(\psi^{\langle n,1\rangle})\in H^{m+1}$, \eqref{Fl} and $L$ depends on $\|\psi^{\langle n, 1\rangle}\|_m$ and $\|\psi^{(n,1)}\|_m$, or equivalently depends on $R$ due to \eqref{semi-eror} and \eqref{eq:psib} by induction.
Hence
\begin{align*}
\|e^{n+1}\|_l&\le e^{C L\varepsilon^{p}\tau}\|e^n\|_l+CM_2(L+M_2^{p})\varepsilon^{p}\tau h^{1+m-l}\\
&\le e^{C L\varepsilon^{p}(n+1)\tau}\|e^0\|_l+CM_2(L+M_2^{p})\varepsilon^{p}\tau h^{1+m-l}\sum\limits_{k=0}^n e^{kC L\varepsilon^{p}\tau}\\
&\le C e^{C LT_0}h^{1+m-l}\|\psi_0\|_{m+1}+\frac{LM_2+M_2^{p+1}}{L}e^{CLT_0}h^{1+m-l}\\
&\le M_3h^{1+m-l},
\end{align*}
where $M_3:=\max\{C_1\|\psi_0\|_{m+1}, C e^{C LT_0}\|\psi_0\|_{m+1}+\frac{LM_2+M_2^{p+1}}{L}e^{CLT_0}\}$ depends on $T_0$, $R$ and $\|\psi_0\|_{m+1}$. The second inequality in \eqref{eq:psib} can be derived by using the triangle inequality and \eqref{semi-eror}:
\[\|I_N(\psi^n)\|_m\leq \|P_N(\psi^{[n]})\|_m+
\|e^n\|_m\le C\|\psi^{[n]}\|_m+M_3 h\le C(1+R)+1,\]
when $h\le h_0:=1/M_3$. Furthermore, it follows from \eqref{eq:psib} that for any $0\le n\le \frac{T_0/\varepsilon^\beta}{\tau}$,
\[\|I_N(\psi^n)\|_{m+1}\le \|P_N(\psi^{[n]})\|_{m+1}+\|e^n\|_{m+1}\le
C\|\psi^{[n]}\|_{m+1}+M_3\le CM_2+M_3,\]
which immediately gives \eqref{reg} by recalling \eqref{ufull}.

Combining \eqref{sp}-\eqref{eq:psib}, we derive for $s\in (1/2, m]$,
\[\|\psi(\cdot, t_n)-I_N(\psi^n)\|_s\le M_1\varepsilon^{p-\beta}\tau^2+M_4h^{1+m-s},\]
where $M_1$ depends on $T_0$, $R$ and $\|\psi\|_{L^\infty([0, T_\varepsilon]; H_{\rm per}^{m+1})}$, and $M_4$ depends on $T_0$, $R$ and $\|\psi_0\|_{m+1}$.
Recalling \eqref{ufull}, we obtain error bounds for $u^n$ and $v^n$ as
\begin{align*}
\|u(\cdot, t_n)-I_N(u^n)\|_s&=\frac{1}{2}\left\|\psi(\cdot, t_n)+\overline{\psi(\cdot, t_n)}-
I_N(\psi^n)-I_N(\overline{\psi^n})\right\|_s\\
&\le \|\psi(\cdot, t_n)-I_N(\psi^n)\|_s\le M_1\varepsilon^{p-\beta}\tau^2+M_4h^{1+m-s},\\
 \|v(\cdot, t_n)-I_N(v^n)\|_{s-1}&=\frac{1}{2}\|\langle \nabla\rangle(\psi(\cdot, t_n)-\overline{\psi(\cdot, t_n)})-\langle \nabla\rangle (I_N(\psi^n)-I_N(\overline{\psi^n}))\|_{s-1}\\
&\le \|\psi(\cdot, t_n)-I_N(\psi^n)\|_s\le M_1 \varepsilon^{p-\beta}\tau^2+M_4h^{1+m-s},
\end{align*}
which shows \eqref{eq:error1} and
the proof for Theorem \ref{thm:eb1} is completed.
\hfill $\square$ \bigskip

\begin{remark}
We  remark here that the same error bounds can be established under the same assumption for the other Strang splitting
\[\psi^{[n+1]}=\mathcal{S}_{\tau}(\psi^{[n]})=\varphi^{\tau/2}_V \circ\varphi^{\tau}_T \circ\varphi^{\tau/2}_V(\psi^{[n]}),\]
and the corresponding full discretization.
Note that
\begin{align*}
\mathcal{S}_{\tau}(\psi_n)&=\varphi^{\tau/2}_V \big[e^{i\tau\langle\nabla\rangle}\psi_n+ \frac{1}{2}\varepsilon^{p}\tau e^{i\tau\langle\nabla\rangle} F(\psi_n)\big]\\
&=e^{i\tau\langle\nabla\rangle}\psi_n+ \frac{1}{2}\varepsilon^{p}\tau e^{i\tau\langle\nabla\rangle} F(\psi_n)+\frac{1}{2}\varepsilon^{p}\tau F\big(e^{i\tau\langle\nabla\rangle}\psi_n+ \frac{1}{2}\varepsilon^{p}\tau e^{i\tau\langle\nabla\rangle} F(\psi_n)\big)\\
&=e^{i\tau\langle\nabla\rangle}\psi_n+ \frac{1}{2}\varepsilon^{p}\tau e^{i\tau\langle\nabla\rangle} F(\psi_n)+\frac{1}{2}\varepsilon^{p}\tau F\big(e^{i\tau\langle\nabla\rangle}\psi_n\big)+E_2,
\end{align*}
where
\[E_2=\frac{1}{4}\varepsilon^{2p}\tau^2\int_0^1
F'\big(e^{i\tau\langle\nabla\rangle}\psi_n+ \frac{\theta}{2}\varepsilon^{p}\tau e^{i\tau\langle\nabla\rangle} F(\psi_n)\big)\big(e^{i\tau\langle\nabla\rangle} F(\psi_n)\big)d\theta.\]
Thus by \eqref{exact}, we get
\be\label{loc2}
\mathcal{S}_{\tau}(\psi_n)-\mathcal{S}_{e, \tau}(\psi_n)
=\varepsilon^{p} e^{i\tau\langle\nabla\rangle}r_3+E_2-\varepsilon^{2p} e^{i\tau\langle\nabla\rangle} r_2-\varepsilon^{3p} e^{i\tau\langle\nabla\rangle}  E_3,
\ee
where
\[r_3=\frac{\tau}{2}\left(F_0(\psi_n)+F_\tau(\psi_n)\right)-\int_0^\tau F_\theta (\psi_n)d\theta=\frac{\tau^3}{2}\int_0^1 \theta (1-\theta)\partial^2_\omega F_\omega (\psi_n)|_{\omega=\theta\tau}d\theta\lesssim \tau^3.\]
It remains to estimate $E_2$. By \eqref{Fp}, we have
\begin{align*}
&\hspace{-4mm}F'\big(e^{i\tau\langle\nabla\rangle}\psi_n+ \frac{\theta}{2}\varepsilon^{p}\tau e^{i\tau\langle\nabla\rangle} F(\psi_n)\big)\big(e^{i\tau\langle\nabla\rangle} F(\psi_n)\big)\\
&=\frac{(p+1)i}{2^{p+1}}\langle \nabla \rangle^{-1}\Big[e^{i\tau\langle\nabla\rangle}\big(\psi_n+ \frac{\theta}{2}\varepsilon^{p}\tau F(\psi_n)\big)+
e^{-i\tau\langle\nabla\rangle}\big(\overline{\psi_n}-
\frac{\theta}{2}\varepsilon^{p}\tau F(\psi_n)\big)\Big]^{p}\\
&\qquad\qquad\qquad\left(e^{i\tau\langle\nabla\rangle} F(\psi_n)-e^{-i\tau\langle\nabla\rangle} F(\psi_n)\right)\\
&=-(p+1)\langle \nabla \rangle^{-1}\Big[\mathrm{Re}\Big(e^{i\tau\langle\nabla\rangle}\big(\psi_n+ \frac{\theta}{2}\varepsilon^{p}\tau F(\psi_n)\big)\Big)\Big]^{p} \sin(\tau\langle\nabla\rangle) F(\psi_n),
\end{align*}
which implies that
\begin{align*}
&\hspace{-4mm}\Big\|F'\big(e^{i\tau\langle\nabla\rangle}\psi_n+ \frac{\theta}{2}\varepsilon^{p}\tau e^{i\tau\langle\nabla\rangle} F(\psi_n)\big)\big(e^{i\tau\langle\nabla\rangle} F(\psi_n)\big)\Big\|_s\\
& \le C\left\|\psi_n+ \frac{\theta}{2}\varepsilon^{p}\tau F(\psi_n)\right\|_s^{p}\;\big\|\sin(\tau\langle\nabla\rangle) F(\psi_n)\big\|_s\\
&\le C\tau \left(\|\psi_n\|_s+\varepsilon^{p}\tau \|F(\psi_n)\|_s\right)^{p}\|F(\psi_n)\|_{s+1}\\
&\le C\tau \left(\|\psi_n\|_s+C\varepsilon^{p}\tau \|\psi_n\|^{p+1}_s\right)^{p}\|\psi_n\|_s^{p+1}\lesssim \tau.
\end{align*}
This suggests that $E_2\lesssim \varepsilon^{2p}\tau^3$, which directly yields that
\[\mathcal{S}_{\tau}(\psi_n)-\mathcal{S}_{e, \tau}(\psi_n)\lesssim \varepsilon^{p}\tau^3.\]
Then the error estimates can be derived by similar and standard arguments.
\end{remark}


\section{Numerical results}
In this section, we present numerical results concerning spatial and temporal accuracy of the TSFP method \eqref{ufull} via \eqref{psifull} for the NKGE \eqref{eq:21}. In our numerical experiments, we take $p=2$, $a=0$ and $b=2\pi$ in \eqref{eq:21} and choose the initial data as
\begin{equation}
u_0(x) = \frac{3}{2}\sin(2x)\quad \mbox{and} \quad u_1(x)  = \frac{5}{1+\cos^2(x)},\quad x \in[0, 2\pi].
\label{eq:long_initial}
\end{equation}
The computation is carried out on a time interval $[0, T_0/\varepsilon^\beta]$ with $0\le \beta\le p=2$ and $T_0>0$ fixed. Here, we study the following three cases with respect to different $\beta$:

(i). Fixed time dynamics up to the time at $O(1)$, i.e., $\beta = 0$;

(ii). Intermediate long-time dynamics up to the time at $O(\varepsilon^{-1})$, i.e., $\beta = 1$;

(ii). Long-time dynamics up to the time at $O(\varepsilon^{-2})$, i.e., $\beta = 2$.

The `exact' solution $u(x, t)$ is obtained numerically by using the TSFP \eqref{psifull}--\eqref{ufull} with a fine mesh size $h_e = \pi/64$ and a very small time step $\tau_e = 10^{-5} $. Denote $u^n_{h, \tau}$ as the numerical solution obtained by the TSFP \eqref{psifull}--\eqref{ufull} with mesh size $h$ and time step $\tau$ at the time $t=t_n$. The errors are denoted as $e(x, t_n) = u(x, t_n) - I_N(u^n_{h, \tau})(x)$. In order to quantify the numerical errors, we measure the $H^1$-norm of $e(\cdot, t_n)$.

\begin{figure}[h!]
\centerline{\psfig{figure=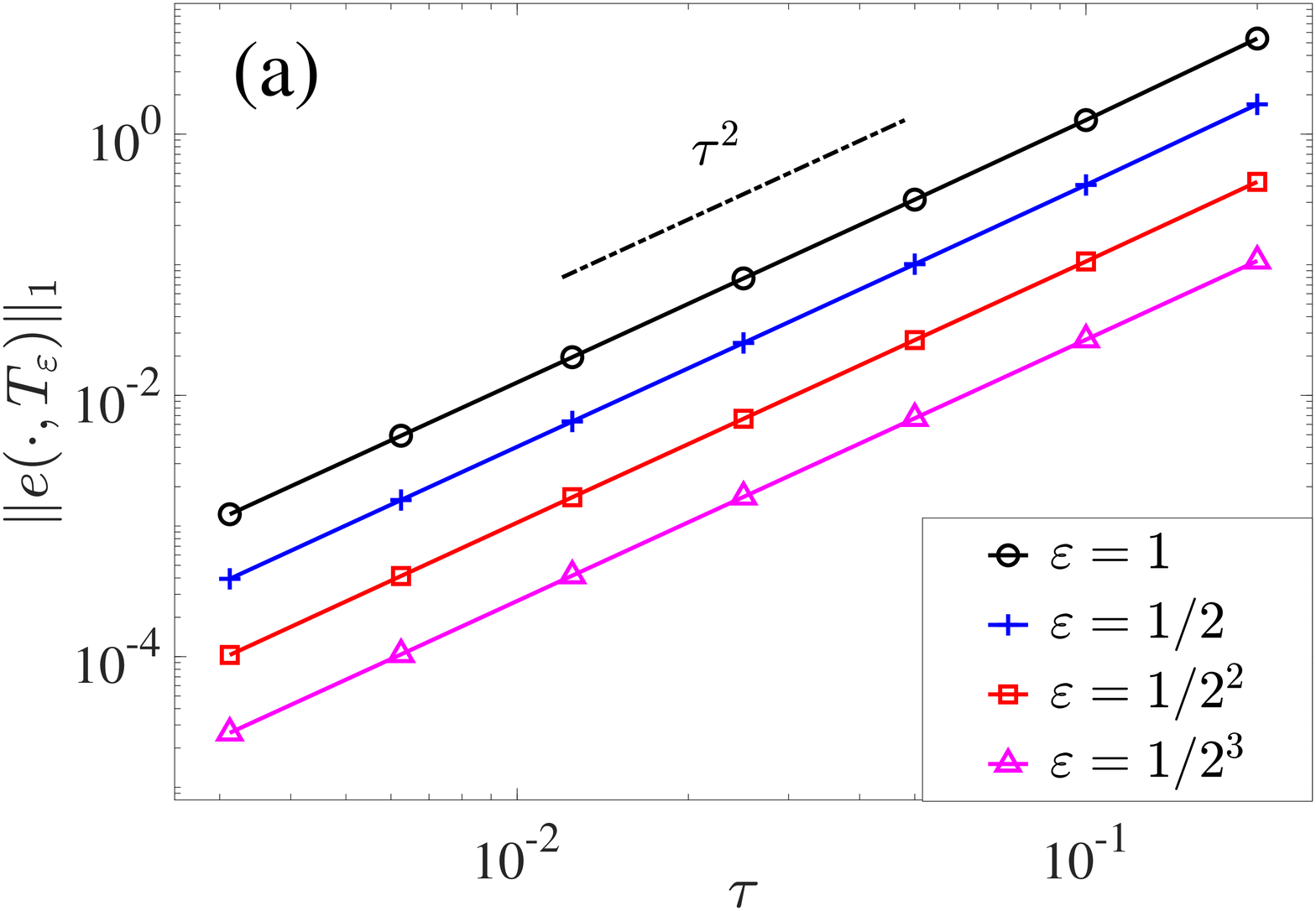,height=5.5cm,width=7cm}
\psfig{figure=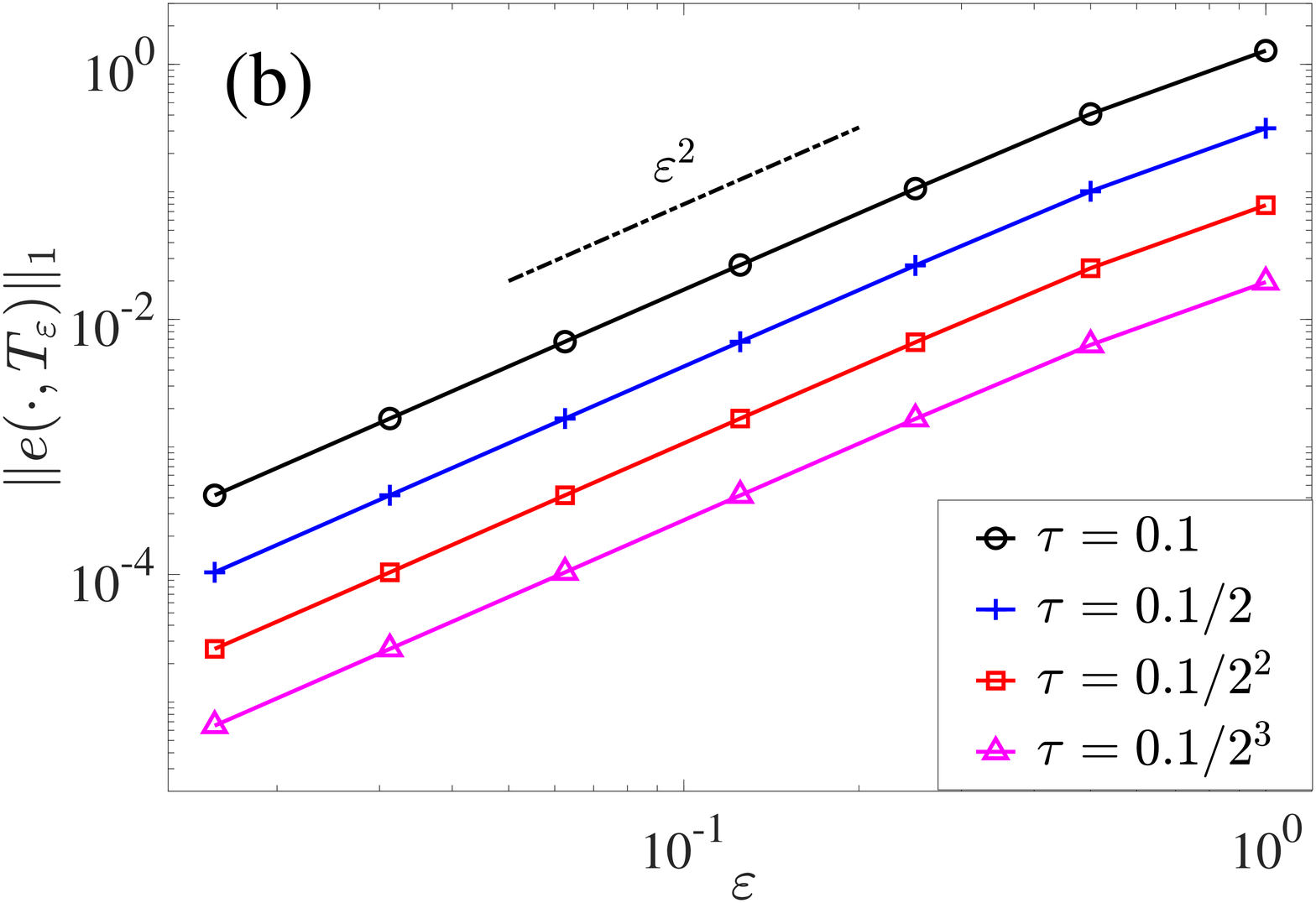,height=5.5cm,width=7cm}}
\caption{Temporal errors of the TSFP \eqref{psifull}--\eqref{ufull} for the NKGE \eqref{eq:21} with $\beta = 0$ for different $\varepsilon$ and $\tau$.}
\label{fig:beta0_t}
\end{figure}

\begin{figure}[h!]
\centerline{\psfig{figure=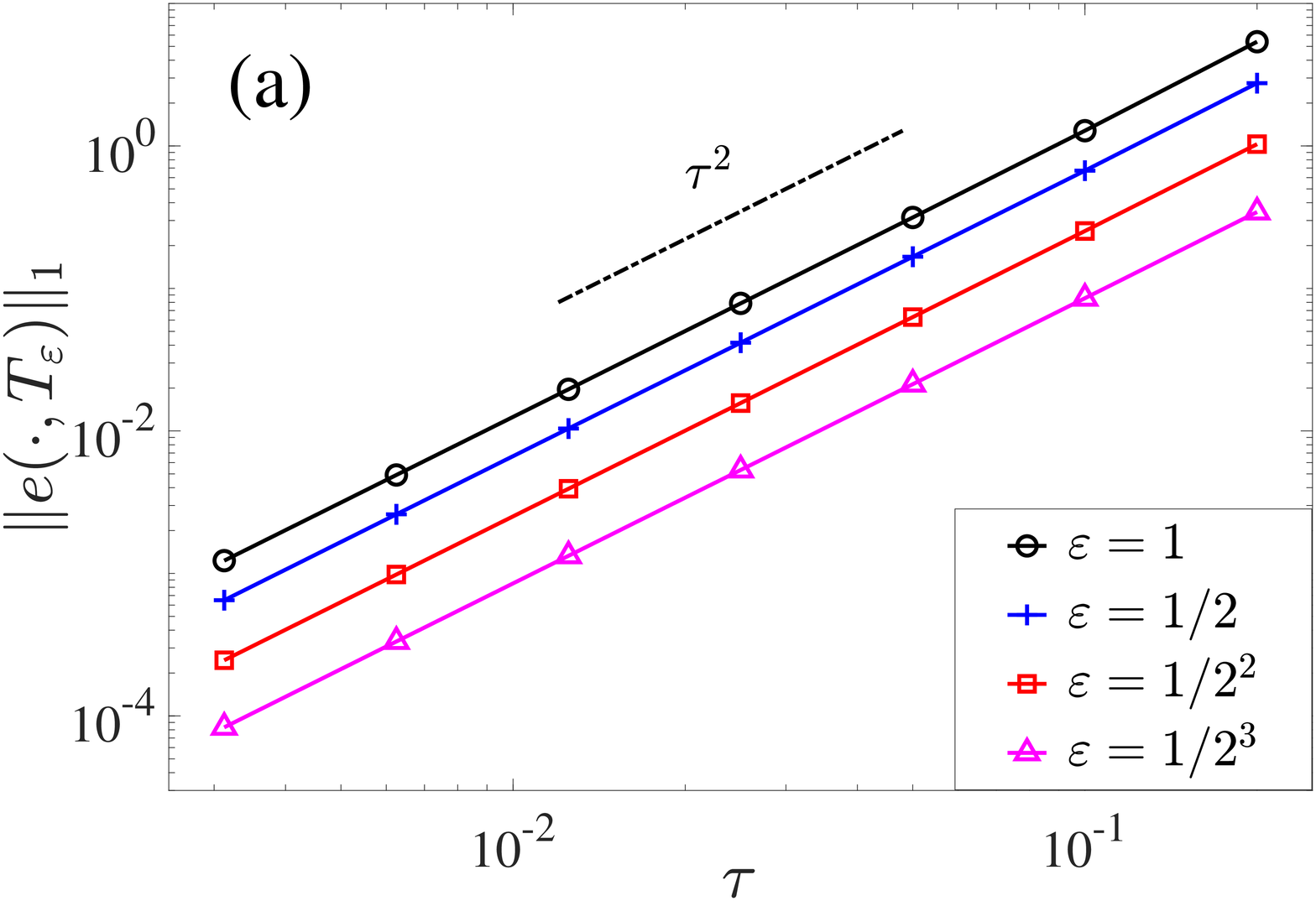,height=5.5cm,width=7cm}
\psfig{figure=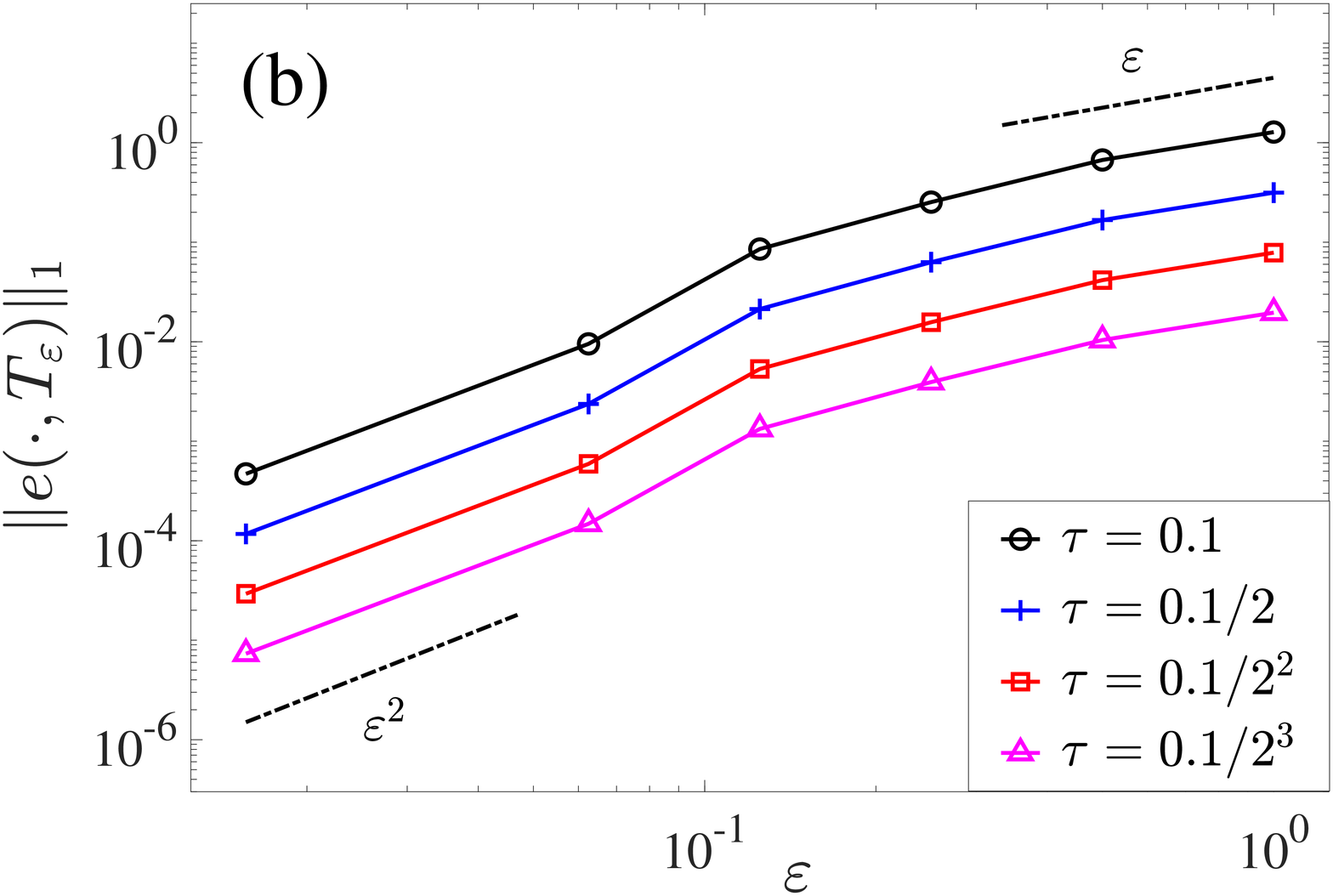,height=5.5cm,width=7cm}}
\caption{Temporal errors of the TSFP \eqref{psifull}--\eqref{ufull} for the NKGE \eqref{eq:21} with $\beta = 1$ for different $\varepsilon$ and $\tau$.}\label{fig:beta1_t}
\end{figure}

\begin{figure}[h!]
\centerline{\psfig{figure=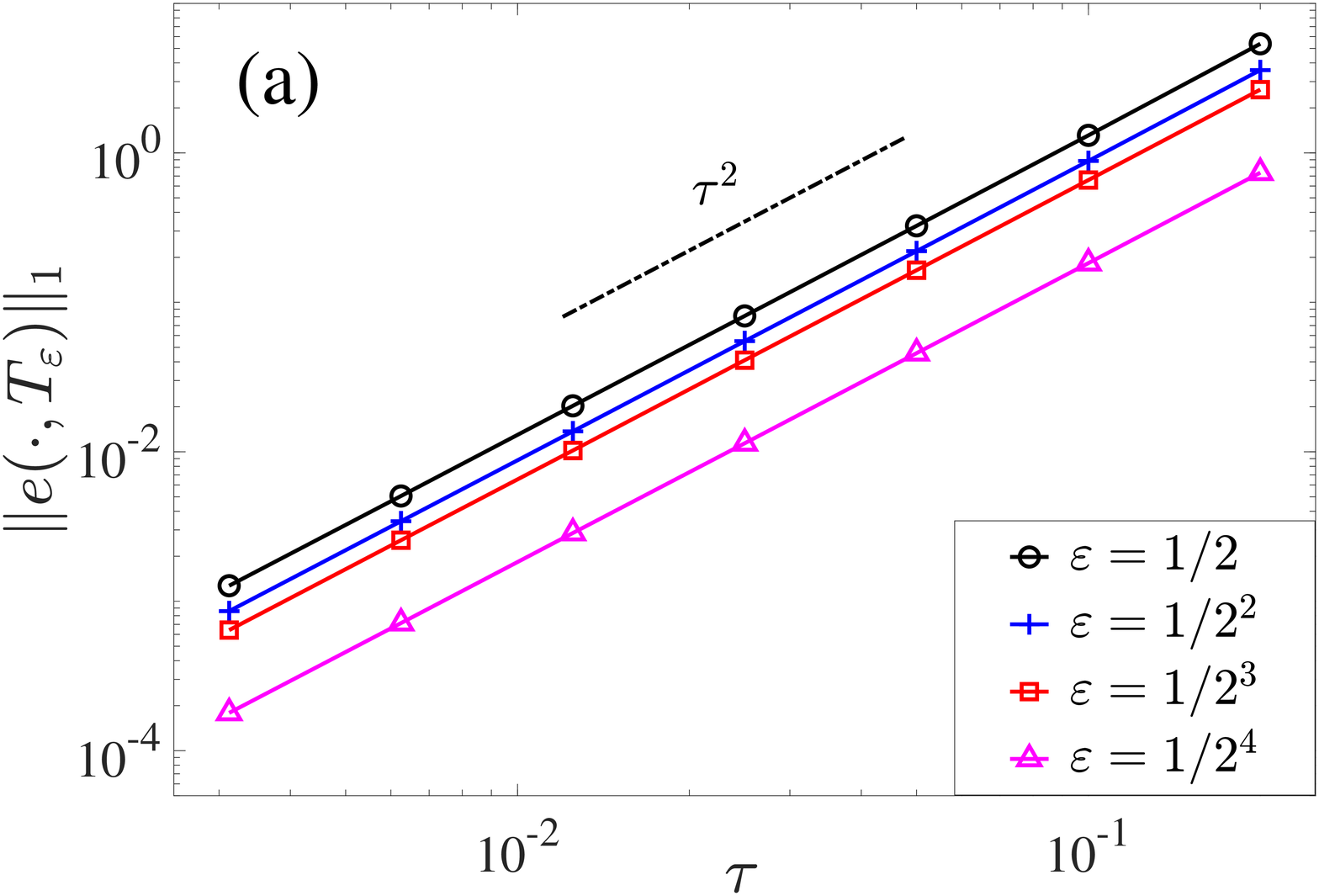,height=5.5cm,width=7cm}
\psfig{figure=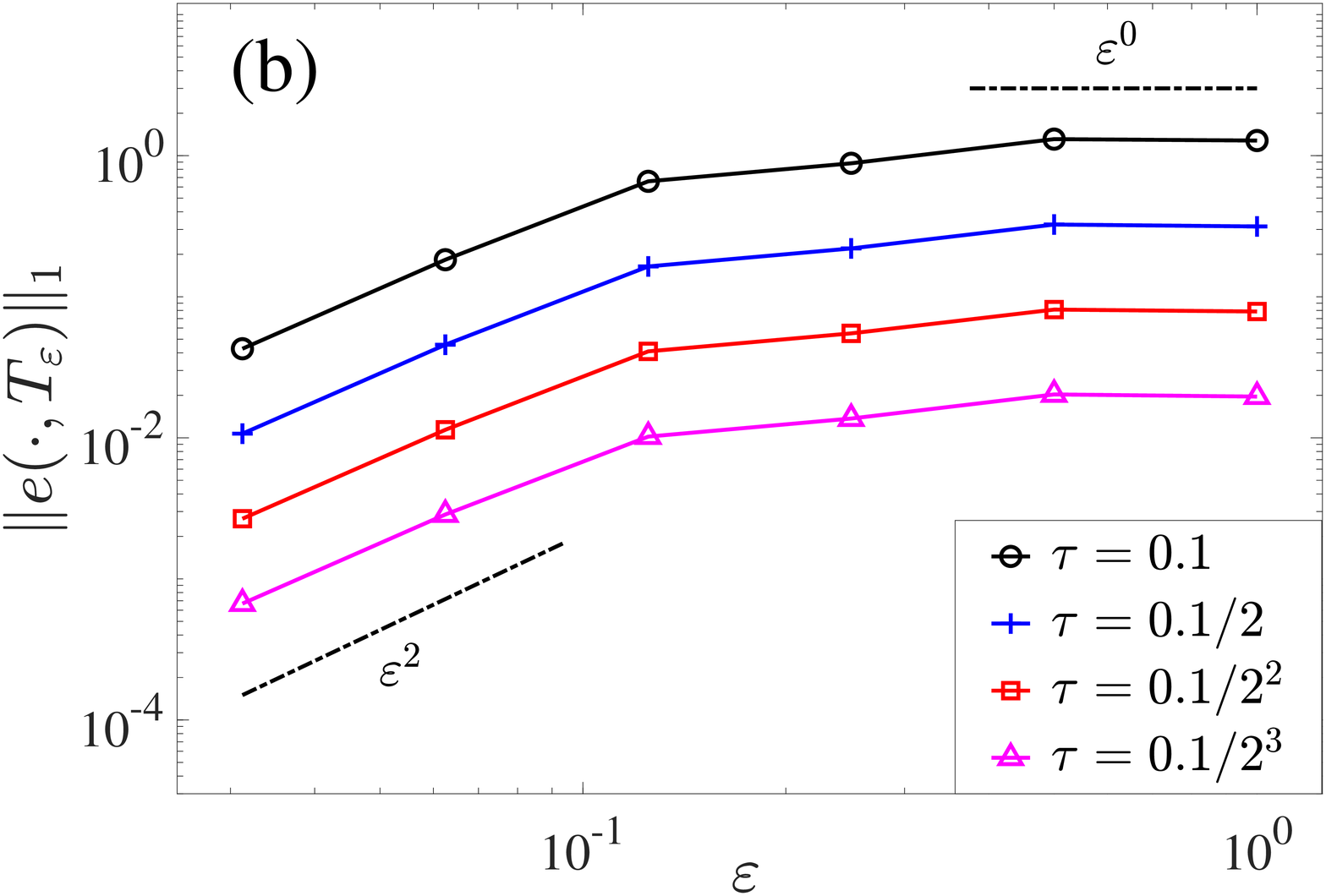,height=5.5cm,width=7cm}}
\caption{Temporal errors of the TSFP \eqref{psifull}--\eqref{ufull} for the NKGE \eqref{eq:21} with $\beta = 2$ for different $\varepsilon$ and $\tau$.}\label{fig:beta2_t}
\end{figure}

The errors are displayed at $T_0 = 1$ with different $\varepsilon$ and $\beta$. For spatial error analysis, we fix the time step as $\tau = 10^{-5}$ such that the temporal errors can be neglected; for temporal error analysis, a very fine mesh size $h = \pi/64$ is chosen such that the spatial errors can be ignored. Table \ref{tab:h} shows the spatial errors under different mesh size $h$ and Figures \ref{fig:beta0_t}--\ref{fig:beta2_t} depict the temporal errors for $\beta = 0$, $\beta = 1$ and $\beta = 2$, respectively.

\begin{table}[h!]
\def\temptablewidth{1\textwidth}
\setlength{\tabcolsep}{6pt}
\caption{Spatial errors of the TSFP \eqref{psifull}--\eqref{ufull} for the NKGE \eqref{eq:21} with \eqref{eq:long_initial} for different $\beta$ and $\varepsilon$.}
\label{tab:h}
{\rule{\temptablewidth}{0.75pt}}
\centering
\begin{tabular*}{\temptablewidth}{@{\extracolsep{\fill}}c|ccccc}
& $\|e(\cdot,T_{\varepsilon})\|_1$ &$h_0 = \pi/4 $ & $h_0/2 $ &$h_0/2^2 $ & $h_0/2^3$  \\
\hline
\multirow{6}{*}{$\beta=0$}
&$\varepsilon_0 = 1$ & 1.12E-1 & 1.22E-3 & 5.03E-6 & 1.54E-12  \\
&$\varepsilon_0 / 2 $ & 8.99E-2 & 6.32E-4 & 2.05E-6 & 1.25E-12  \\
&$\varepsilon_0 / 2^2 $ & 9.04E-2 & 4.67E-4 & 1.95E-6 & 1.19E-12  \\
&$\varepsilon_0 / 2^3 $ & 8.85E-2 & 4.47E-4 & 1.93E-6 & 1.18E-12  \\
&$\varepsilon_0 / 2^4 $ & 8.82E-2 & 4.47E-4 & 1.93E-6 & 1.19E-12 \\
&$\varepsilon_0 / 2^5 $ & 8.81E-2 & 4.48E-4 & 1.93E-6 & 1.18E-12 \\
\hline
\multirow{6}{*}{$\beta=1$}
&$\varepsilon_0 = 1$ & 1.12E-1 & 1.22E-3 & 5.03E-6 & 1.54E-12  \\
&$\varepsilon_0 / 2 $ & 2.14E-1 & 2.10E-3 & 1.58E-6 & 5.72E-13 \\
&$\varepsilon_0 / 2^2 $ & 1.08E-1 & 2.36E-3 & 7.09E-7 & 1.24E-12 \\
&$\varepsilon_0 / 2^3 $ & 4.47E-2 & 9.27E-4 & 7.72E-7 & 1.52E-13 \\
&$\varepsilon_0 / 2^4 $ & 1.14E-1 & 8.11E-4 & 7.13E-7 & 7.97E-13 \\
&$\varepsilon_0 / 2^5 $ & 7.29E-2 & 1.24E-3 & 9.83E-7 & 1.26E-12 \\
\hline
\multirow{6}{*}{$\beta=2$}
&$\varepsilon_0 = 1$ & 1.12E-1 & 1.22E-3 & 5.03E-6 & 1.54E-12  \\
&$\varepsilon_0 / 2 $ & 5.22E-1 & 6.58E-3 & 5.81E-7 & 1.16E-12 \\
&$\varepsilon_0 / 2^2 $ & 5.79E-1 & 1.52E-3 & 1.82E-6 & 1.20E-12 \\
&$\varepsilon_0 / 2^3 $ & 5.82E-1 & 1.03E-3 & 6.05E-7 & 9.90E-13 \\
&$\varepsilon_0 / 2^4 $ & 9.17E-1 & 1.68E-3 & 6.69E-7 & 4.78E-12 \\
&$\varepsilon_0 / 2^5 $ & 7.67E-1 & 1.79E-3 & 3.52E-7 & 1.22E-11 \\
\end{tabular*}
{\rule{\temptablewidth}{0.75pt}}
\end{table}

From Table \ref{tab:h} and Figures \ref{fig:beta0_t}--\ref{fig:beta2_t}, we can draw the following observations:

(1) The TSFP method converges uniformly for $0 < \varepsilon \leq 1$ in space with exponential convergence rate (cf. each row in Table \ref{tab:h}).

(2) For any fixed $\varepsilon = \varepsilon_0 > 0$, the TSFP method \eqref{psifull}--\eqref{ufull} is second-order accurate in time (cf. each line in Figures \ref{fig:beta0_t}(a)--\ref{fig:beta2_t}(a)). When $\beta=0$, the temporal error behaves like $O(\varepsilon^2\tau^2)$ (cf. Figure \ref{fig:beta0_t}(b)), which agrees with the theoretical result in Theorem \ref{thm:eb1}. Figure \ref{fig:beta1_t}(b) and Figure \ref{fig:beta2_t}(b) show that the temporal error is at $O(\varepsilon\tau^2)$ and $O(\tau^2)$ for $\beta=1$ and $\beta=2$, respectively.

(3) Our numerical results confirm the uniform error bounds in Theorem \ref{thm:eb1}.

\newcolumntype{L}[1]{>{\raggedright\arraybackslash}p{#1}}
\newcolumntype{C}[1]{>{\centering\arraybackslash}p{#1}}
\newcolumntype{R}[1]{>{\raggedleft\arraybackslash}p{#1}}

\begin{table}[h!]
\caption{Comparison of the properties of the TSFP method \eqref{psifull}--\eqref{ufull} for the NKGE
\eqref{eq:WNE} at different long-time dynamics regimes with $T_0$ and
$\tau_0$ are fixed and independent of $0<\varepsilon\le 1$.}
\label{tab:comp_long}
\vspace*{-10pt}
\begin{center}
\begin{tabular}{C{2.5cm}|C{2.5cm}|C{2.5cm}|C{2.5cm}|C{2.5cm}}
\toprule
  & $O(1)$-time with $\beta=0$    & intermediate long-time with $0 < \beta \leq \frac{p}{3}$
  &long-time with $\frac{p}{3} < \beta < p$ & super long-time with  $\beta = p$ \\
\hline
final time $T_{\varepsilon}=\frac{T_0}{\varepsilon^{\beta}}$ &\vspace{0.2mm} $O(1)$ & long-time $O(\varepsilon^{-\beta})$ &longer-time  $O(\varepsilon^{-\beta})$ & longest-time  $O(\varepsilon^{-p})$ \\
\hline
largest time step size $\tau_{\varepsilon} = \tau_0 \varepsilon^{\frac{\beta - p}{2}}$ & largest time step size at $O(\varepsilon^{-p/2})$ & larger time \quad  step size at $O(\varepsilon^{\frac{\beta-p}{2}})$ & large time \quad  step size at $O(\varepsilon^{\frac{\beta-p}{2}})$ & \hspace{2cm} $O(1)$ \\
  \hline
total time steps $N_{\varepsilon} = \frac{T_{\varepsilon}}{\tau_{\varepsilon}} = \frac{T_0}{\tau_0}\varepsilon^{\frac{p-3\beta}{2}}$  &\hspace{2cm} $O(1)$ &\hspace{2cm} $O(1)$ &\hspace{2cm} $O(\varepsilon^{\frac{p-3\beta}{2}})$ & \hspace{2cm} $O(\varepsilon^{-p})$ \\ \hline
total computational cost &\hspace{2cm} $O(N\ln N)$ &\hspace{2cm} $O(N\ln N)$ &\hspace{2cm} $O(\varepsilon^{\frac{p-3\beta}{2}}N\ln N)$ & \hspace{2cm} $O(\varepsilon^{-p}N\ln N)$   \\
\hline
\vspace{0.15mm}spatial  error &uniform spectral &uniform spectral &uniform spectral &uniform spectral \\
 \hline
temporal error in term of $\tau_0$ & uniform second-order  & uniform second-order & uniform second-order & uniform second-order \\
\bottomrule
\end{tabular}
\end{center}
\end{table}

\smallskip

We remark here that, when $0<\varepsilon\ll1$, our numerical results suggest a better error bound as (cf. left-bottom parts in Figures
\ref{fig:beta1_t}(b)--\ref{fig:beta2_t}(b))
\[\|u(\cdot, t_n) - I_N(u^n)\|_s\lesssim h^{1+m-s} + \varepsilon^p \tau^2,\quad 0 \leq n \leq \frac{T_0/\varepsilon^\beta}{\tau}.\]
We emphasized this improved convergence when $\varepsilon\ll 1$ is missing for the EWI-FP method presented in \cite{FY}, which shows the superiority of the TSFP method.

For convenience of readers, Table \ref{tab:comp_long} summarizes
the properties of the TSFP method \eqref{psifull}--\eqref{ufull} for the NKGE
\eqref{eq:WNE} at different long-time dynamics regimes.

\section{Extension to an oscillatory complex NKGE in the whole space}
In this section, we begin with a complex NKGE in the whole space,
re-scale it into an oscillatory complex NKGE, compare properties of the
NKGE under different scalings and extend the TSFP method and its error bounds
to the oscillatory complex NKGE.

\subsection{Comparisons of the complex NKGE under different scalings}

Consider the following complex NKGE with a power-type nonlinearity in the whole space $\mathbb{R}^d$ ($d=1,2,3$) as
\be\label{eq:WNEws}
\left\{
\begin{aligned}
&\partial_{tt}u({\bx}, t)-\Delta u({\bx}, t)+ u({\bx}, t)+ \varepsilon^{p}|u({\bx}, t)|^{p}u({\bx}, t)=0,\quad{\bx} \in \mathbb{R}^d,\quad t > 0,\\
&u({\bx}, 0)=u_0({\bx})=O(1),\quad\partial_t u({\bx},0)=u_1({\bx})=O(1),\quad{\bx} \in \mathbb{R}^d.
\end{aligned}\right.
\ee
Here, $u:= u(\bx, t)$ is a complex-valued scalar field, and the initial datum $u_0({\bx})$ and $u_1({\bx})$ are two given complex-valued functions which are independent of the parameter $\varepsilon$. Again formally, the amplitude of the solution $u$ is at $O(1)$. The local/global well-posedness of the Cauchy problem \eqref{eq:WNEws} and scattering properties have been extensively studied in a considerable literature \cite{GV,KT,K,Lan,MS1972,MS2018,Nak1999,OTT}. Particularly, under appropriate assumptions on $p$, $d$, $\varepsilon$ and the initial conditions, the solutions of \eqref{eq:WNEws} are global \cite{CN} and scatter as $|t|\rightarrow \infty$ for small initial values (low energy scattering) \cite{K,MS2018}, or for all initial values (asymptotic completeness) \cite{MS1972,Nak1999}.
 In addition, under proper regularity of the solution, the complex NKGE \eqref{eq:WNEws} is time symmetric or time reversible and conserves the energy \cite{BD,BFY,DXZ} as
\begin{align}
E_1(t) &:=E_1(u(\cdot,t))= \int_{\mathbb{R}^d} \left[ |\partial_t u (\bx, t)|^2 + |\nabla u(\bx, t)|^2 + |u(\bx, t)|^2 +\frac{2\varepsilon^{p}}{p+2} |u(\bx, t)|^{p+2}  \right] d\bx\nn\\
&\equiv \int_{\mathbb{R}^d} \left[ |u_1(\bx)|^2 + |\nabla u_0(\bx)|^2 + |u_0(\bx)|^2 +\frac{2\varepsilon^{p}}{p+2} |u_0(\bx)|^{p+2}  \right] d\bx \label{eq:Energy_uws}\\
&= E_1(0) = O(1), \qquad t \geq 0.\nn
\end{align}
Plugging the plane wave solution $u(\bx, t) = Ae^{i(\boldsymbol \xi \cdot {\bx}-\omega_1 t)}$ (with $A$ the amplitude, $\boldsymbol{\xi}$ the spatial wave number  and $\omega_1:=\omega_1(\boldsymbol{\xi})$ the time frequency) into the complex NKGE \eqref{eq:WNEws}, we
get the dispersion relation:
\be\label{drws}
\omega_1=\omega_1({\boldsymbol{\xi}}) =
\pm\sqrt{1+|{\boldsymbol{\xi}}|^2+\varepsilon^{p} A^p}=O(1), \quad \varepsilon\in (0, 1], \quad \mathrm{for\,\,\,fixed}\,\,\,\boldsymbol{\xi}\in\mathbb{R}^d,
\ee
which immediately implies the group velocity
\be\label{gvws}
\boldsymbol{v}_1:=\boldsymbol{v}_1(\xi)=\nabla \omega_1(\boldsymbol \xi)=\pm\frac{\boldsymbol\xi}{\sqrt{1+|\boldsymbol\xi|^2+ \varepsilon^p  A^p}}= O(1).
\ee
Thus the solution of the complex NKGE \eqref{eq:WNEws} propagates waves with
amplitude at $O(1)$, wavelength in space and time at $O(1)$ and wave velocity at $O(1)$.

By introducing $w({\bx}, t)=\varepsilon u({\bx}, t)$, we can reformulate the complex NKGE \eqref{eq:WNEws} with weak nonlinearity (and initial data with amplitude at $O(1)$)  into the following complex NKGE with small initial data (and $O(1)$ nonlinearity):
\begin{equation}
\label{eq:SIEws}
\left\{
\begin{aligned}
&\partial_{tt} w({\bx}, t)-\Delta w({\bx}, t)+ w({\bx}, t)+ |w({\bx}, t)|^{p}w({\bx}, t)=0, \quad \bx \in \mathbb{R}^d,\quad t > 0, \\
&w({\bx}, 0) = \varepsilon u_0({\bx})=O(\varepsilon),\quad \partial_t w({\bx}, 0) = \varepsilon u_1({\bx})=O(\varepsilon),\quad {\bx} \in \mathbb{R}^d.
\end{aligned}\right.
\end{equation}
Noticing that the amplitude of the initial data in
\eqref{eq:SIEws} is at $O(\varepsilon)$, formally we can get the amplitude of the solution $w$ of \eqref{eq:SIEws} is at $O(\varepsilon)$, too. Similarly, the complex NKGE \eqref{eq:SIEws} is time symmetric or time reversible and conserves the energy \cite{BD,BFY,DXZ} as
\begin{align*}
E_2(t) &:=E_2(w(\cdot,t))= \int_{\mathbb{R}^d} \big[ |\partial_t w (\bx, t)|^2 + |\nabla w(\bx, t)|^2 + |w(\bx, t)|^2 +\frac{2}{p+2} |w(\bx, t)|^{p+2}  \big] d \bx\\
&\equiv \varepsilon^2 \int_{\mathbb{R}^d} \left[ | u_1(\bx)|^2 + |\nabla u_0(\bx)|^2 + | u_0(\bx)|^2 +\frac{2\varepsilon^{p}}{p+2}|u_0(\bx)|^{p+2}  \right]  d\bx\\
&=E_2(0)=\varepsilon^2 E_1(0)= O(\varepsilon^2), \qquad t\ge0.
\end{align*}
In addition, plugging the plane wave solution $w(\bx, t) = \varepsilon Ae^{i(\boldsymbol \xi \cdot {\bx}-\omega_1 t)}$  into the complex NKGE \eqref{eq:SIEws}, we get the same dispersion relation \eqref{drws} and
the same group velocity \eqref{gvws} of the complex NKGE \eqref{eq:SIEws},
i.e., the complex NKGEs  \eqref{eq:SIEws} and \eqref{eq:WNEws} share
the  same dispersion relation \eqref{drws} and
the same group velocity \eqref{gvws}.
Again, the solution of the complex NKGE \eqref{eq:SIEws} propagates waves with
amplitude at $O(\varepsilon)$, wavelength in space and time at $O(1)$ and wave velocity  at $O(1)$.

Introducing a re-scale in time
\be\label{rstime}
t=\frac{s}{\varepsilon^{\beta}}\Leftrightarrow s=\varepsilon^\beta t, \qquad \nu(\bx,s)=u(\bx,t),
\ee
with $0 < \beta \leq p$ fixed, we can re-formulate the complex NKGE \eqref{eq:WNEws}
into the following oscillatory complex NKGE
\begin{equation}
\left\{
\begin{split}
&\partial_{ss} \nu(\bx, s)+\frac{1}{\varepsilon^{2\beta}}(-\Delta+1) \nu(\bx, s) + \varepsilon^{p-2\beta}|\nu(\bx, s)|^{p}\nu(\bx, s)=0,\,\,\,\, \bx \in \mathbb{R}^d,\,\,\,\, s > 0,\\
&\nu(\bx, 0) = u_0(\bx)=O(1),\quad \partial_s \nu(\bx, 0)={\varepsilon}^{-\beta}u_1(\bx)=O({\varepsilon}^{-\beta}),\quad \bx \in \mathbb{R}^d.
\end{split}\right.
\label{eq:HOEws}
\end{equation}
Formally, the amplitude of the solution $\nu$ of the oscillatory complex NKGE \eqref{eq:HOEws} is at $O(1)$. Again, the oscillatory complex NKGE \eqref{eq:HOEws} is time symmetric or time reversible and conserves the energy \cite{BD,BFY,DXZ} as
\begin{align}
E_3(s) &:=E_3(\nu(\cdot,s))= \int_{\mathbb{R}^d} \left[ |\partial_s \nu|^2 +\frac{1}{\varepsilon^{2\beta}} \left(|\nabla \nu|^2 + |\nu|^2\right) +\frac{2\varepsilon^{p-2\beta}}{p+2} |\nu|^{p+2}  \right] d\bx\nn\\
&\equiv \frac{1}{\varepsilon^{2\beta}}\int_{\mathbb{R}^d} \left[ |u_1(\bx)|^2 + |\nabla u_0(\bx)|^2 + |u_0(\bx)|^2 +\frac{2\varepsilon^{p}}{p+2} |u_0(\bx)|^{p+2}  \right] d\bx \label{eq:Energy_uws3}\\
&= E_3(0) =\varepsilon^{-2\beta}E_1(0)= O(\varepsilon^{-2\beta}), \qquad s \geq 0.\nn
\end{align}
Again, plugging the plane wave solution $\nu(\bx, s) = Ae^{i(\boldsymbol \xi\cdot {\bx}-\omega_2 s)}$ (with $A$ the amplitude, $\boldsymbol\xi$ the spatial wave number  and $\omega_2:=\omega_2(\boldsymbol\xi)$ the time frequency) into the oscillatory complex NKGE \eqref{eq:HOEws}, we
get the dispersion relation:
\be\label{drhoe}
\omega_2=\omega_2(\boldsymbol\xi) =
\pm\frac{1}{\varepsilon^\beta}\sqrt{1+|\boldsymbol\xi|^2+ \varepsilon^{p} A^p}=O(\varepsilon^{-\beta}), \qquad \boldsymbol\xi\in\mathbb{R}^d,
\ee
which immediately implies the group velocity
\be\label{gvhoe}
\boldsymbol v_2:=\boldsymbol v_2(\boldsymbol\xi)=\nabla \omega_2(\boldsymbol \xi)=\pm\frac{\boldsymbol\xi}{\varepsilon^\beta \sqrt{1+|\boldsymbol\xi|^2+\varepsilon^{p} A^p}}= O(\varepsilon^{-\beta}).
\ee
Thus the solution of the oscillatory complex NKGE \eqref{eq:HOEws}
propagates waves with amplitude at $O(1)$, wavelength in space and time at $O(1)$ and $O(\varepsilon^\beta)$, respectively, and wave velocity at
$O(\varepsilon^{-\beta})$.

\begin{remark}
We remark here that the above scalings of the complex NKGE are different
from the following complex NKGE in the nonrelativistic regime, which has been widely used and studied in the literature \cite{BCZ,BD,BMS2004,BZ,BFS,MN,Sch1979}:
\be\label{eq:WNEwsnrr}
\left\{
\begin{aligned}
&\partial_{tt}u({\bx}, t)-\frac{1}{\varepsilon^2}\Delta u({\bx}, t)+ \frac{1}{\varepsilon^4}u({\bx}, t)+\frac{1}{\varepsilon^2}|u({\bx}, t)|^2u({\bx}, t)=0,\quad t > 0,\\
&u({\bx}, 0)=u_0({\bx})=O(1),\quad\partial_t u({\bx},0)=\varepsilon^{-2}u_1({\bx})=O(\varepsilon^{-2}),\quad{\bx} \in \mathbb{R}^d.
\end{aligned}\right.
\ee
The above complex NKGE conserves the energy \cite{BD,BFY,DXZ} as
\begin{align}
E_4(t) &:=E_4(u(\cdot,t))=
\int_{\mathbb{R}^d} \Big[ |\partial_t u (\bx, t)|^2 + \frac{|\nabla u(\bx, t)|^2}{\varepsilon^2} +\frac{1}{\varepsilon^4} |u(\bx, t)|^2 +\frac{1}{2\varepsilon^2} |u(\bx, t)|^4  \Big] d\bx\nn\\
&\equiv \frac{1}{\varepsilon^4}\int_{\mathbb{R}^d} \left[ |u_1(\bx)|^2 + \varepsilon^2|\nabla u_0(\bx)|^2 + |u_0(\bx)|^2 +\frac{\varepsilon^2}{2} |u_0(\bx)|^4  \right] d\bx \label{eq:Energy_uwsrr}\\
&= E_5(0) = O(\varepsilon^{-4}), \qquad t \geq 0.\nn
\end{align}
Plugging the plane wave solution $u(\bx, t) = Ae^{i(\boldsymbol\xi \cdot {\bx}-\omega_3 t)}$ into the complex NKGE \eqref{eq:WNEwsnrr}, we
get the dispersion relation:
\be\label{drwsrr1}
\omega_3=\omega_3(\boldsymbol\xi) =
\pm\frac{1}{\varepsilon^2}\sqrt{1+\varepsilon^2|\boldsymbol\xi|^2+\varepsilon^2 A^2}=O(\varepsilon^{-2}), \qquad\boldsymbol\xi\in\mathbb{R}^d,
\ee
which immediately implies the group velocity
\be\label{gvwsrr1}
\boldsymbol v_3:=\boldsymbol v_3(\boldsymbol\xi)=\nabla \omega_3(\boldsymbol \xi)=\pm\frac{\boldsymbol\xi}{\sqrt{1+\varepsilon^2|\boldsymbol\xi|^2+\varepsilon^2 A^2}}= O(1).
\ee
Thus the solution of the complex NKGE \eqref{eq:WNEwsnrr} propagates waves with amplitude at $O(1)$, wavelength in space and time at $O(1)$ and $O(\varepsilon^2)$, respectively, and wave velocity at $O(1)$.
\end{remark}

For convenience of readers, Table \ref{tab:comp1} shows the properties of the complex NKGE under different scalings.

\begin{table}[http]
\centering
\renewcommand\arraystretch{1.2}
  \caption{Comparison of the complex NKGE under different scalings.}\label{tab:comp1}
  \vspace*{-10pt}
\def\temptablewidth{1\textwidth}
{\rule{\temptablewidth}{0.75pt}}
\begin{tabular*}{\temptablewidth}{@{\extracolsep{\fill}}c|p{2.3cm}<{\centering}|p{2.3cm}<{\centering}|p{2.3cm}<{\centering}|p{2.3cm}<{\centering}}
  &\eqref{eq:WNEws}    &\eqref{eq:SIEws}
  &\eqref{eq:HOEws} &\eqref{eq:WNEwsnrr}  \\[0.25em]
\hline
amplitude &$O(1)$ &$O(\varepsilon)$ &$O(1)$  &$O(1)$  \\ \hline
wavelength in space &$O(1)$ &$O(1)$ &$O(1)$  &$O(1)$  \\  \hline
wavelength in time &$O(1)$ &$O(1)$ &$O(\varepsilon^{\beta})$ &$O(\varepsilon^{2})$   \\ \hline
wave velocity &$O(1)$ &$O(1)$ &$O(\varepsilon^{-\beta})$ &$O(1)$   \\  \hline
energy &$O(1)$ &$O(\varepsilon^2)$ &$O(\varepsilon^{-2\beta})$ &$O(\varepsilon^{-4})$
\end{tabular*}
{\rule{\temptablewidth}{0.75pt}}
\end{table}

\subsection{The TSFP method for the complex NKGE \eqref{eq:HOEws} and main results}
Similar to those in the literature, we truncate
the oscillatory complex NKGE \eqref{eq:HOEws} in 1D
onto a bounded interval $\Omega = (a, b)$ with periodic boundary conditions as
\be\left\{
\begin{split}
&\partial_{ss}\nu(x, s)+\frac{1}{\varepsilon^{2\beta}}(-\partial_{xx}+1)\nu(x, s) + \varepsilon^{p-2\beta} |\nu(x, s)|^{p} \nu(x, s)= 0,\,\,\, s > 0,\\
&\nu(a,t)=\nu(b,t),\quad \partial_x \nu(a,t)=\partial_x \nu(b,t), \qquad t\ge0,\\
&\nu(x, 0)=u_0(x),\quad \partial_s \nu(x, 0)=\varepsilon^{-\beta} u_1(x),\quad x \in \overline{\Omega}=[a, b].
\end{split}\right.
\label{eq:HOE_w1}
\ee

Denote $q(x, s) =\partial_s \nu(x, s)$, by taking $k=\varepsilon^{\beta}\tau$ and assuming $u_0$ and $u_1$ to be real-valued in \eqref{eq:HOE_w1},  the TSFP discretization can be similarly obtained via \eqref{psifull}. Under the following reasonable assumptions on the exact solution $\nu$ of the oscillatory NKGE \eqref{eq:HOE_w1}
\[
{\rm(B)}
\begin{split}
&\nu \in  \  L^\infty\left([0, T_0]; H^{m+1}_{\rm per}\right), \qquad
\partial_s \nu\in L^\infty\left([0, T_0]; H^{m}_{\rm per}\right),\\
&\|\nu\|_{L^{\infty}\left([0, T_0]; H^{m+1}_{\rm per}\right)} \lesssim 1,\qquad\,\,\,\, \|\partial_s \nu\|_{L^{\infty}\left([0, T_0]; H^{m}_{\rm per} \right)} \lesssim \frac{1}{\varepsilon^\beta},
\end{split}
\]
with $m\ge 1$, we can establish the following error bounds of the TSFP method for the oscillatory complex NKGE \eqref{eq:HOE_w1} (the proof is omitted here for brevity).

\begin{theorem}
\label{thm:eb1_HOE}
Let $\nu^n$, $q^n$ be the numerical approximation obtained from the TSFP method. Under the assumption (B), there exist $h_0 > 0$ and $k_0 > 0$ sufficiently small and independent of $\varepsilon$ such that, for any $0 < \varepsilon \leq 1$, when $0 < h \leq h_0$ and $0 < k \leq k_0 \varepsilon^{\frac{3\beta-p}{2}}$,  we have the error estimates for $l \in (1/2, m]$
\[\|\nu(\cdot, s_n) - I_N(\nu^n)\|_l+\varepsilon^\beta\|\partial_s \nu(\cdot, s_n) - I_N(q^n)\|_{l-1} \lesssim h^{1+m-l} + \varepsilon^{p-3\beta}k^2,\quad 0 \leq n \leq \frac{T_0}{k}.
\]
\end{theorem}

\begin{remark}
From Theorem \ref{thm:eb1_HOE}, we clearly see that the TSFP is uniformly second-order accurate in the weakly oscillatory case, i.e., $0\le \beta \leq \frac{p}{3}$. Furthermore, large time step size at $k\sim \varepsilon^{\frac{3\beta-p}{2}}$ is allowed in practical computation when $0\le \beta<\frac{p}{3}$. While for $\beta\in (\frac{p}{3}, p]$, the TSFP method fails to be uniformly convergent and tiny time step is required as $k\lesssim \varepsilon^{\frac{3\beta-p}{2}}$.
\end{remark}

\subsection{Numerical results}
In order to verify the error bounds in Theorem \ref{thm:eb1_HOE}, we take $d=1$ and $p=3$ in \eqref{eq:HOEws} and the initial data
\be
u_0(x)=(2+i)e^{-x^2/2},\quad u_1(x)=\, \sech(x^2),\quad x \in \mathbb{R}.
\label{eq:HOE_initial}
\ee
The problem is solved on a bounded interval $\Omega_{\varepsilon} = [-8 - \varepsilon^{-\beta}, 8 + \varepsilon^{-\beta}]$ since the wave velocity is
at $O(\varepsilon^{-\beta})$, which is large enough to guarantee that the periodic boundary condition does not introduce a significant truncation error relative to the original problem. The `exact' solution $v(x, s)$ is obtained numerically by using the TSFP method with a fine mesh size $h_e = 1/16$ and a very small time step $k_e = {10}^{-5}$. We also measure the $H^1$-norm and the errors are displayed at $T_0 = 1$ with different $\varepsilon$ and $\beta$. For the oscillatory complex NKGE \eqref{eq:HOEws}, we study the following three cases:

Case I. Weakly oscillatory regime, i.e., $\beta = 1$;

Case II. Intermediate oscillatory regime, i.e., $\beta = 2$;

Case III. Highly oscillatory regime, i.e., $\beta = 3$.

For spatial error analysis, we fix the time step as $k ={10}^{-5}$ such that the temporal errors can be neglected; for temporal error analysis, a very fine mesh size $h = 1/16$ is chosen such that the spatial error can be ignored. Table \ref{tab:HOE_h} shows the spatial errors under different mesh size for these three cases and Tables \ref{tab:HOE_t1}-\ref{tab:HOE_t3} depict the temporal errors for $\beta=1, 2, 3$, respectively. In order to quantify the error, we introduce
\[ e_{\infty} (t) := \max_{0 < \varepsilon \leq 1} \{\|e(\cdot, t)\|_1\}.\]

\begin{table}[h!]
\def\temptablewidth{1\textwidth}
\setlength{\tabcolsep}{6pt}
\caption{Spatial errors of the TSFP method  for the oscillatory complex NKGE \eqref{eq:HOE_w1} with \eqref{eq:HOE_initial} for different $\beta$ and $\varepsilon$.}
\label{tab:HOE_h}
{\rule{\temptablewidth}{0.75pt}}
\centering
\begin{tabular*}{\temptablewidth}{@{\extracolsep{\fill}}c|ccccc}
& $\|e(\cdot,1)\|_1$ &$h_0 = 1 $ & $h_0/2 $ &$h_0/2^2 $ & $h_0/2^3$  \\
\hline
\multirow{4}{*}{$\beta=1$}
&$\varepsilon_0 =1/2 $ & 1.57E-1 & 2.40E-3 & 5.82E-6 & 1.77E-9  \\
&$\varepsilon_0 / 2^1 $ & 9.52E-2 & 2.91E-3 & 8.47E-6 & 2.59E-10  \\
&$\varepsilon_0 / 2^2 $ & 5.85E-2 & 3.31E-3 & 1.11E-5 & 3.32E-10  \\
&$\varepsilon_0 / 2^3 $ & 1.03E-1 & 1.63E-3 & 1.18E-5 & 4.10E-10  \\
\hline
\multirow{4}{*}{$\beta=2$}
&$\varepsilon_0 =1/2 $ & 2.01E-1 & 3.24E-3 & 1.02E-5 & 1.52E-9  \\
&$\varepsilon_0 / 2^1 $ & 1.11E-1 & 1.64E-3 & 1.19E-5 & 4.34E-10  \\
&$\varepsilon_0 / 2^2 $ & 1.28E-1 & 3.57E-3 & 1.55E-5 & 1.61E-10  \\
&$\varepsilon_0 / 2^3 $ & 1.18E-1 & 3.81E-3 & 1.34E-5 & 2.21E-10  \\
\hline
\multirow{4}{*}{$\beta=3$}
&$\varepsilon_0 =1/2 $ & 1.91E-1 & 3.90E-3 & 1.40E-5 & 6.95E-9  \\
&$\varepsilon_0 / 2^1 $ & 1.55E-1 & 3.43E-3 & 1.58E-5 & 1.66E-10  \\
&$\varepsilon_0 / 2^2 $ & 1.30E-1 & 5.79E-3 & 5.94E-6 & 4.32E-10  \\
&$\varepsilon_0 / 2^3 $ & 1.25E-1 & 5.15E-3 & 1.59E-5 & 5.95E-10  \\
\end{tabular*}
{\rule{\temptablewidth}{0.75pt}}
\end{table}

\begin{table}[h!]
\def\temptablewidth{1\textwidth}
\setlength{\tabcolsep}{6pt}
\caption{Temporal errors of the TSFP method  for the oscillatory complex NKGE \eqref{eq:HOE_w1} with \eqref{eq:HOE_initial} and $\beta=1$.}
\label{tab:HOE_t1}
{\rule{\temptablewidth}{0.75pt}}
\centering
\begin{tabular*}{\temptablewidth}{@{\extracolsep{\fill}}cccccccc}
$\|e(\cdot, 1)\|_1$ &$k_0 = 0.1 $ & $k_0/2 $ &$k_0/2^2 $ & $k_0/2^3$ & $k_0/2^4$  & $k_0/2^5$ & $k_0/2^6$ \\
\hline
$\varepsilon_0 = 1$ & 2.82E-1 & 6.71E-2 & 1.66E-2 & 4.13E-3 & 1.03E-3 & 2.58E-4 & 6.45E-5 \\
order & - & 2.07 & 2.02 & 2.01 & 2.00 & 2.00 & 2.00 \\
\hline
$\varepsilon_0 / 2 $ & 1.15E-1 & 2.77E-2 &6.85E-3 & 1.71E-3 & 4.27E-3 & 1.07E-4 & 2.67E-5 \\
order & -  & 2.05 & 2.02 & 2.00 & 2.00 & 2.00 & 2.00 \\
\hline
$\varepsilon_0 / 2^2 $ & 4.20E-2 & 9.45E-3 & 2.31E-3 & 5.75E-4 & 1.43E-4 & 3.58E-5 & 8.96E-6 \\
order & -  & 2.15 & 2.03 & 2.01 & 2.01 & 2.00  & 2.00 \\
\hline
$\varepsilon_0 / 2^3 $ & 4.91E-2 & 6.43E-3 & 1.46E-3 & 3.57E-4 & 8.89E-5 & 2.22E-5 & 5.54E-6 \\
order & -  & 2.93 & 2.14 & 2.03 & 2.01 & 2.00 & 2.00 \\
\hline
$\varepsilon_0 / 2^4 $ & 2.29E-2 & 8.02E-3 & 1.01E-3 & 2.29E-4 & 5.60E-5 & 1.39E-5 & 3.48E-6 \\
order & -  & 1.51 & 2.99 & 2.14 & 2.03 & 2.01 & 2.00  \\
\hline
$\varepsilon_0 / 2^5 $ & 8.77E-3 & 3.48E-3 & 1.21E-3 & 1.51E-4 & 3.43E-5 & 8.40E-6 & 2.09E-6 \\
order & -  & 1.33 & 1.52 & 3.00 & 2.14 & 2.03 & 2.01\\
\hline
$\varepsilon_0 / 2^6 $ & 9.87E-4 & 1.25E-3 & 4.88E-4 & 1.70E-4 & 2.10E-5 & 4.78E-6 & 1.17E-6 \\
order & -  & -0.34 & 1.36 & 1.52  & 3.02 & 2.14 & 2.03\\
\hline\hline
$e_{\infty}(t=1)$ & 2.82E-1 & 6.71E-2 & 1.66E-2 & 4.13E-3 & 1.03E-3 & 2.58E-4 & 6.45E-5 \\
order & - & 2.07 & 2.02 & 2.01 & 2.00 & 2.00 & 2.00 \\
\end{tabular*}
{\rule{\temptablewidth}{0.75pt}}
\end{table}

\begin{table}[h!]
\def\temptablewidth{1\textwidth}
\setlength{\tabcolsep}{6pt}
\caption{Temporal errors of the TSFP method  for the oscillatory complex NKGE \eqref{eq:HOE_w1} with \eqref{eq:HOE_initial} and $\beta=2$.}
\label{tab:HOE_t2}
{\rule{\temptablewidth}{0.75pt}}
\centering
\begin{tabular*}{\temptablewidth}{@{\extracolsep{\fill}}ccccccc}
$\|e(\cdot, 1)\|_1$ &$k_0 = 0.1 $ & $k_0/2 $ &$k_0/2^2 $ & $k_0/2^3$ & $k_0/2^4$  & $k_0/2^5$ \\\hline
$\varepsilon_0 = 1$ & \bf{2.82E-1} & 6.71E-2 & 1.66E-2 & 4.13E-3 & 1.03E-3 & 2.58E-4 \\
order & \bf{-} & 2.07 & 2.02 & 2.01 & 2.00 & 2.00 \\
\hline
$\varepsilon_0 / 4^{1/3} $ & 5.15E-1 & \bf{1.14E-1} & 2.77E-2 & 6.89E-3 & 1.72E-3 & 4.30E-4 \\
order & -  & \bf{2.18} & 2.04 & 2.01 & 2.00 & 2.00 \\
\hline
$\varepsilon_0 / 4^{2/3} $ & 1.52 & 2.20E-1 & \bf{5.08E-2} & 1.25E-2 & 3.11E-3 & 7.76E-4 \\
order & -  & 2.79 & \bf{2.11} & 2.02 & 2.01 & 2.00 \\
\hline
$\varepsilon_0 / 4 $ & 1.40 & 6.80E-1 & 8.95E-2 & \bf{2.03E-2} & 4.96E-3 & 1.23E-3 \\
order & -  & 1.04 & 2.93 & \bf{2.14} & 2.03 & 2.01 \\
\hline
$\varepsilon_0 / 4^{4/3} $ & 9.33E-1 & 6.94E-1 & 3.18E-1 & 3.88E-2 & \bf{8.11E-3} & 1.96E-3 \\
order & -  & 0.43 & 1.13 & 3.03 & \bf{2.26} & 2.05  \\
\hline
$\varepsilon_0 / 4^{5/3} $ & 3.10E-1 & 2.48E-1 & 2.85E-1 & 1.19E-1 & 2.07E-2 & \bf{3.45E-3} \\
order & -  & 0.32 & -0.20 & 1.26 & 2.52 & \bf{2.58} \\
\end{tabular*}
{\rule{\temptablewidth}{0.75pt}}
\end{table}

\begin{table}[h!]
\def\temptablewidth{1\textwidth}
\setlength{\tabcolsep}{6pt}
\caption{Temporal errors of the TSFP method for the oscillatory complex NKGE \eqref{eq:HOE_w1} with \eqref{eq:HOE_initial} and $\beta=3$.}
\label{tab:HOE_t3}
{\rule{\temptablewidth}{0.75pt}}
\centering
\begin{tabular*}{\temptablewidth}{@{\extracolsep{\fill}}ccccccc}
$\|e(\cdot, 1)\|_1$ &$k_0 = 0.1 $ & $k_0/4 $&$k_0/4^2 $&$k_0/4^3$&$k_0/4^4$ & $k_0/4^5$\\
\hline
$\varepsilon_0 = 1$ & \bf{2.82E-1} & 1.66E-2 & 1.03E-3 & 6.45E-5 & 4.03E-6 & 2.50E-7 \\
order & \bf{-} & 2.04 & 2.01 & 2.00 & 2.00 & 2.01 \\
\hline
$\varepsilon_0 / 4^{1/3} $ & 3.82 & \bf{1.34E-1} & 8.23E-3 & 5.14E-4 & 3.21E-5 & 1.99E-6 \\
order & -  & \bf{2.42} & 2.01 & 2.00 & 2.00 & 2.01 \\
\hline
$\varepsilon_0 / 4^{2/3} $ & 8.46 & 6.37E-1 & \bf{3.35E-2} & 2.08E-3 & 1.30E-4 & 8.03E-6 \\
order & -  & 1.87 & \bf{2.12} & 2.00 & 2.00 & 2.01 \\
\hline
$\varepsilon_0 / 4 $ & 4.08 & 1.95 & 1.22E-1 & \bf{6.76E-3} & 4.20E-4 & 2.60E-5 \\
order & -  & 0.53 & 2.00 & \bf{2.09} & 2.00 & 2.01 \\
\hline
$\varepsilon_0 / 4^{4 /3}$ & 1.39 & 1.15 & 5.40E-1 & 2.61E-2 & \bf{1.45E-3} & 8.92E-5 \\
order & -  & 0.14 & 0.55 & 2.19 & \bf{2.08} & 2.01 \\
\hline
$\varepsilon_0 / 4^{5/3} $ & 4.26E-1 & 3.59E-1 & 2.98E-1 & 1.39E-1 & 6.17E-3 & \bf{3.40E-4} \\
order & -  & 0.12 & 0.13 & 0.55 & 2.25 & \bf{2.09} \\
\end{tabular*}
{\rule{\temptablewidth}{0.75pt}}
\end{table}

From Tables \ref{tab:HOE_h}-\ref{tab:HOE_t2} and additional similar results not shown here for brevity, we can draw the following observations for  the TSFP method:

(1) The TSFP method is uniformly and spectrally accurate in space for any $0\le\beta \leq p$ (cf. Table \ref{tab:HOE_h}).

(2) When $\beta = 1$, the TSFP method converges quadratically in time, which is uniformly for $0<\varepsilon\le 1$ (cf. last row in Table \ref{tab:HOE_t1}).
While for cases $\beta = 2$ and $\beta=3$, second-order convergence can only be observed when $k\lesssim \varepsilon^{3/2}$ and $k\lesssim \varepsilon^3$, respectively (cf. the upper triangle above the main diagonal in Tables \ref{tab:HOE_t2}-\ref{tab:HOE_t3}). This agrees with the analytical result in Theorem \ref{thm:eb1_HOE}.

(3) For $\beta=2$ and $\beta=3$, when $0<\varepsilon \ll1$ and $k\lesssim \varepsilon^{\frac{3\beta-p}{2}}$, Tables \ref{tab:HOE_t2}--\ref{tab:HOE_t3} suggest the following improved error bound
\[
\|v(\cdot, s_n) - I_N(v^n)\|_l\lesssim h^{1+m-l} + \varepsilon^{p-2\beta}k^2,\quad 0 \leq n \leq T_0/k.
\]

Again, for convenience of readers, Table \ref{tab:comp2} summarizes the properties of the TSFP method for the oscillatory NKGE \eqref{eq:HOEws} at different parameter regimes.

\begin{table}[ht]
\caption{Comparison of the properties of the TSFP method
for the oscillatory NKGE \eqref{eq:HOE_w1} at different parameter
regimes, while $T_0$ and
$k_0$ are fixed and independent of $0<\varepsilon\le 1$.}
\label{tab:comp2}
\begin{center}
\begin{tabular}{m{1.7cm}|m{1.3cm}|m{1.6cm}|m{1.3cm}|m{1.8cm}|m{1.8cm}|m{1.8cm}}
\toprule
  & $\beta = 0$    & $ 0 < \beta < \frac{p}{3}$
  &$\beta = \frac{p}{3}$ &$\frac{p}{3} < \beta < \frac{p}{2}$ & $\beta = \frac{p}{2}$ & $\frac{p}{2} < \beta \leq p$ \\ [3pt]
\hline
nonlinearity strength $\varepsilon^{p-2\beta}$ & weakest $O(\varepsilon^p)$ & weaker $O(\varepsilon^{p-2\beta})$ & weak $O(\varepsilon^{p/3})$ & weak $O(\varepsilon^{p-2\beta})$ & $O(1)$ & strong $O(\varepsilon^{p-2\beta})$ \\
\hline
 \multicolumn{1}{m{1.7cm}|}{time step size $k_{\varepsilon} = k_0 \varepsilon^{\frac{3\beta - p}{2}}$} &  \multicolumn{1}{m{1.3cm}|}{larger $O(\varepsilon^{-p/2})$} & large $O(\varepsilon^{\frac{3\beta-p}{2}})$ & $O(1)$ & small $O(\varepsilon^{\frac{3\beta-p}{2}})$ & smaller  $O(\varepsilon^{p/4})$ & smallest $O(\varepsilon^{\frac{3\beta-p}{2}})$   \\
  \hline
total time steps $N_{\varepsilon} = \frac{T_0}{k_{\varepsilon}} = \frac{T_0}{k_0}\varepsilon^{\frac{p-3\beta}{2}}$  &$O(1)$ &$O(1)$ & $O(1)$ & many steps at $O(\varepsilon^{\frac{p-3\beta}{2}})$ & many steps at $O(\varepsilon^{-\frac{p}{4}})$ &  many steps at $O(\varepsilon^{\frac{p-3\beta}{2}})$   \\ \hline
total cost &\small $O(N\ln N)$ &\small  $O(N\ln N)$ &\small  $O(N\ln N)$ &\small  $O\left(\frac{N\ln N}{\varepsilon^{\frac{3\beta-p}{2}}}\right)$ &\small  $O\left(\varepsilon^{\frac{p}{4}}N\ln N\right)$& \small  $O\left(\frac{N\ln N}{\varepsilon^{\frac{3\beta-p}{2}}}\right)$   \\
 \hline
\tabincell{c}{spatial \\ error} & \tabincell{c}{uniform \\ spectral}
& \tabincell{c}{uniform \\ spectral} & \tabincell{c}{uniform \\ spectral}
& \tabincell{c}{uniform \\ spectral} & \tabincell{c}{uniform \\ spectral}
& \tabincell{c}{uniform \\ spectral}
 \\
\hline
temporal error & uniform & uniform & uniform &   non-uniform &   non-uniform &   non-uniform \\
\bottomrule
\end{tabular}
\end{center}
\end{table}

\section{Conclusion}
An efficient and accurate time-splitting Fourier pseudospectral (TSFP) method was proposed and analyzed for the long-time dynamics of the nonlinear Klein--Gordon equation (NKGE) with weak nonlinearity or small initial data. Uniform error bounds of the TSFP method were established up to the time at $O(\varepsilon^{-p})$ with $0<\varepsilon\le 1$ a dimensionless parameter used to characterize the nonlinearity strength. Numerical results were reported to confirm our error bounds in the long-time regime. Extension of the method and its error bounds to an oscillatory complex NKGE in the whole space was discussed.

\black

\bibliographystyle{amsplain}

\end{document}